\documentstyle[epsf]{article}
\begin{document}           
\title{On the groupoid of transformations of rigid structures on
  surfaces
\footnote{First version: {\tt November 9, 1997}.
This version: {\tt June 30, 1999}.This preprint is available electronically at 
          \tt  http://www-fourier.ujf-grenoble.fr/\~{ }funar }}

\author{
\begin{tabular}{cc}
 Louis Funar &  R{\u{a}}zvan Gelca\\
\small \em Institut Fourier BP 74 
&\small \em Department of Mathematics, 2072 East Hall\\
\small \em University of Grenoble I &\small \em University of Michigan\\
\small \em 38402 Saint-Martin-d'H\`eres cedex, France
&\small \em Ann Arbor MI 48109, USA\\
\small \em e-mail: {\tt funar@fourier.ujf-grenoble.fr}
& \small \em e-mail: {\tt rgelca@math.lsa.umich.edu } \\
\end{tabular}
}

\date{}
\maketitle

{\abstract We prove that
the 2-groupoid of transformations of rigid structures on surfaces 
 has 
a finite presentation, establishing a 
result first conjectured by Moore and Seiberg. We also show that a
finite dimensional, unitary, cyclic topological quantum field theory
gives rise to a representation of this $2$-groupoid.        }
\vspace{10mm}



\tableofcontents

\section{Introduction}

Three dimensional topological quantum field theories
(TQFT's) give rise to representations of the mapping
class groups of closed surfaces. TQFT's with corners
give rise to representations of a related object,
the $2$-groupoid of transformations of rigid structures.
Rigid structures (also called DAP-decompositions in 
\cite{Fro, Gel2, Wa}) are decompositions of 
surfaces into disks, annuli and pairs of pants, together
with additional information for keeping track of 
twistings.

The $2$-groupoid of transformations of rigid structures appeared
for the first time in the works of 
physicists studying $2$-dimensional conformal field theories.
Specifically G. Moore and N. Seiberg (see \cite{MoSe}) 
worked with this groupoid and conjectured a presentation of it. 
In an unpublished preprint \cite{Wa}, K. Walker sketched some ideas for
the proof that the  presentation given by
 Moore and Seiberg  is complete.  As Walker pointed out,
the Moore-Seiberg equations represent compatibility conditions
that the basic data of a TQFT with corners must satisfy.
Based on  Walker's point of view, several TQFT's with 
corners have been constructed so far \cite{Fro, Gel2, Gel}. 
In a TQFT with corners the quantum invariants of 
$3$-manifolds are computed from an initial
amount of information, by making use of the axioms. Of course this initial
amount of information, called basic data, must satisfy the above mentioned
compatibility conditions. Hence the necessity for a rigorous 
proof of the fact that the Moore-Seiberg equations are complete.
This is the purpose of the present paper. In addition to this
we also show how a maximal TQFT (i.e. one that has an underlying
theory with corners), gives rise in a canonical way to a representation
of the $2$-groupoid.

The idea of the proof is to apply the Cerf theoretic techniques 
used by A. Hatcher
and W. Thurston  \cite{HaTh} for obtaining a presentation of the 
mapping class group of a surface. 
Let us mention that an explicit presentation for the mapping class group was 
derived  afterwards by B. Wajnryb (see \cite{Wj,BW}) and a more 
symmetric (but infinite) presentation was given by S. Gervais
(\cite{Ger}). 

The proof given below is done in three
steps. First we exhibit a presentation for the groupoid
of transformations of markings (maximal collections of 
non-isotopic simple closed curves in the interior of 
a surface). Then we explain how this presentation produces
a presentation of the groupoid of overmarkings (collection
of curves cutting a surface in disks, annuli and pairs of
pants). Finally, we use Walker's approach to solve the
case of rigid structures. The last part of the paper describes the
construction of the canonical representations of the $2$-groupoid 
that arise from  TQFT's for which 
 the mapping class group acts in a homogeneous manner. 
We mention that 
our initial result for the case of the 
complex associated to  cut systems  was obtained independently 
 in \cite{HPS}, using the same methods. After this paper appeared in
 preprint form we learned about the work of Bakalov and Kirillov Jr.
 \cite{BK} in which a different proof for the main result is given.
Although their proof is still based on the Hatcher-Thurston ideas, the
authors avoided the direct use of Cerf's theory and use instead
results from \cite{Harer} about cut systems. 

Before proceeding with the details of the paper,  we want to make
some remarks. The $2$-groupoid of transformations
of rigid structures is  a universal
object containing the mapping class 
groups of all surfaces. One can think of
it as playing the role of  the tower of the mapping class groups
of surfaces, a notion suggested by A.Grothendieck in his 
``Esquisse d'un programme''. A more precise connection with 
Grothendieck's program is the relationship between the 
Teichmuller tower of (orbifold) fundamental groupoids of the 
moduli spaces of punctured curves 
and our 2-groupoid (which should be a quotient of the former). 
The basepoints in the moduli spaces 
are chosen in simply connected neigborhoods of infinity, corresponding to the 
maximal semistable degeneracy curves. 
In the context  of topological quantum field theory, 
instead of considering a series of representations 
of mapping class groups, we consider the 
representation of this single but more complicated algebraic  object. 
Notice that this groupoid as a natural central extension related to
those of the mapping class groups (see \cite{MR}). The
representations arising from the most interesting TQFTs are rather
representations of the latter extension. 

{\bf Acknowledgements.} Part of this work when the 
first author visited  Columbia University, whose support and hospitality are
kindly acknowledged. We are thankful to Ch. Frohman, S.Gervais,
T. Kitano, P. Lochak, V. Sergiescu, L. Schneps, J. Stasheff, A.Voronov
 and the referee 
 for useful comments and discussions.

\section{2-groupoids}
\subsection{Algebraic definitions}

\newtheorem{defi}{Definition}[section]
A 1-groupoid is by definition a category whose morphisms are
isomorphisms. We extend this to an object having both the features 
of a 2-category and of a groupoid, and which we will call a 
2-groupoid. 
\newtheorem{defi2}[defi]{Definition}
\begin{defi2}
A 2-groupoid $C$ is a category with the following properties: 
\begin{enumerate}
\item The collection of objects ${\cal O}(C)$ is a category itself, 
which is a 1-groupoid with an associative composition law  denoted by
$\otimes$, which gives ${\cal O}(C)$ the structure of a (strict) 
tensor category. This means that the objects in ${\cal O}(C)$ are 
the homomorphism sets $Hom^0(u,v)$ of some other category $C^0$ having
an associative multiplication. The composition 
$Hom^0(u,v)\times Hom^0(v,w)\to Hom^0(u,w)$ is our tensor structure 
$\otimes$ at the level of ${\cal O}(C)$.
%Set $s, t$ be the source and target maps 
%${\cal O}(C)\longrightarrow F$ into the collection of final objects.
\item On the collection of morphisms one has a
  composition $\circ$ which makes it into a groupoid,
  and a tensor multiplication  
\[\otimes: Hom(X,X') \otimes Hom(Y,Y')\longrightarrow Hom(X\otimes Y,X'\otimes
Y'), \]
induced by $\otimes$ on ${\cal O}(C)$ and compatible with the composition.
Notice that the Hom on $C$ is like a 2-Hom of $C^0$.  
\end{enumerate}
\end{defi2}
%\newtheorem{exe}{Example}[section]
%\begin{exe}
%A simple example is the 2-category $C$ 
%obtained from a category $C_0$ by letting
%${\cal O}(C)=Hom(C_0)$ and 
%$Hom(C)=Isom(Hom(C_0))$.
%\end{exe}
The example we had in mind when considering this definition was that
of the 2-groupoid of  transformations of rigid structures on surfaces.
Recall that a
 DAP-decomposition of a surface $\Sigma $ is a decomposition of the surface
into a finite number of elementary surfaces: disks, annuli, and 
pairs of pants, determined by a collection
 of disjoint simple closed curves in the interior of
$\Sigma $. 
A rigid structure consists of a DAP-decomposition 
together with the  following additional structure:
\begin{enumerate}
\item an ordering of the elementary surfaces;
\item for each elementary surface $\Sigma _0$  a numbering of its
 boundary components, by 1 if $\Sigma _0$ is a disk, 1 and 2 
if $\Sigma _0$ is an annulus, and 1, 2 and 3 if $\Sigma _0$ is a pair
of pants;
\item a parametrization of each boundary component $C$ of 
$\Sigma _0$ by $S^1=\{z;  |z|=1\}$ 
(the parameterization being compatible with the
orientation of $\Sigma _0$ under the convention ``first out'') such that 
the parameterizations coming from two neighboring elementary surfaces are
one the complex conjugate of the other;
\item fixed disjoint embedded arcs
in $\Sigma _0$ joining $e^{i\epsilon}$ (where  $\epsilon >0$ is small)
on the $j$-th boundary component to
$e^{-i\epsilon}$ on the $j+1$-st (modulo the number of boundary
components of $\Sigma _0$) (these arcs are called seams).
\item an ordering of elementary surfaces in the DAP-decomposition
according to topological type.
\end{enumerate}
Rigid structures are considered up to isotopy.
In this setting the category $C^0$ is given by circles (with some
additional structure), and rigid structures on surfaces are 
homomorphisms (in $C^0$) between their boundary. The exterior 
composition on $C^0$ is given by the disjoint union. 
These are related both to the PROPs formalism and to that of the 
modular operads. 
\begin{defi}
The (full) duality groupoid ${\cal D}$ -- also called 
the groupoid of transformations of rigid structures on  
surfaces --   consists of:
\begin{enumerate}
\item A collection of objects 
 $(\Sigma, r)$, which are the rigid two dimensional cobordisms. 
Here
 $\Sigma$ is a surface, with boundary $\partial \Sigma$ endowed
  with a fixed splitting $\partial \Sigma=\partial_+\Sigma \cup 
\partial_-\Sigma$ of the boundary components and two labelings 
of each connected component in $\partial_+\Sigma$ and $\partial_-\Sigma$, and 
$r$ is the rigid structure on $\Sigma$. 
\item The collection of morphisms between two given  objects 
$(\Sigma,r)$ and $(\Sigma',r')$ is the set of all pairs 
$\lambda=(\varphi, c)$, where $\varphi:\Sigma\longrightarrow \Sigma'$ 
is a homeomorphism preserving the boundary splitting (and thus 
$\Sigma=\Sigma'$) and $c:\varphi(r)\longrightarrow r'$ is a change of
the rigid structure. We  factor out by the following 
equivalence relation:
\begin{enumerate}
\item $(\varphi, c)\sim (\varphi', c)$ if $\varphi$ and 
$\varphi'$ are isotopic; 
\item $(\varphi, c)\sim (\varphi', c')$ if
  $c'=c\varphi_*(\varphi'_*)^{-1}$, where $\varphi_*$ is the map induced
  by the homeomorphism $\varphi$ at the level of rigid structures.
\end{enumerate} 
\item
% Notice that all  morphisms 
%$Hom((\Sigma,r),(\Sigma',r')$ are invertible (the setup for a groupoid)
%and have a
The  natural composition  of morphisms, and a 
 tensor product operation such that
 \begin{enumerate}
 \item 
At the level of the objects the (incomplete) tensor 
product is given by:
$(\Sigma,r)\otimes (\Sigma',r')= (\Sigma\otimes \Sigma', r\otimes
  r')$, where $\Sigma\otimes \Sigma'$ is the boundary connected sum 
of $\Sigma$ and $\Sigma'$, identifying the last $k$ connected
components of $\partial_-\Sigma$ with the first $k$ connected components
of $\partial_+\Sigma'$ (here one should
think that the boundary components are
ordered lexicographically according to the
ordering of basic surfaces and to that 
of boundary components within one elementary surface). 
One labels in a canonical way  the connected 
components of  
$\partial_+(\Sigma\otimes \Sigma')$ (which is the union of
of the unglued components of $\partial_+\Sigma$ and  $\partial_+\Sigma'$), and 
likewise  the  components of 
$\partial_-(\Sigma\otimes \Sigma')$. The number $k$
is a parameter of the tensor product. Further 
$r\otimes r'$ is the natural rigid structure induced by the gluing. 
 \item On the  level of morphisms, the tensor product  induces maps:
\[ Hom((\Sigma,r),(\Sigma',r')) \otimes Hom((\tilde\Sigma,\tilde r),
(\tilde\Sigma',\tilde r'))\longrightarrow 
 Hom((\Sigma,r)\otimes\tilde\Sigma,\tilde r) ,
(\Sigma',r')\otimes (\tilde\Sigma',\tilde r' )) \]
 \end{enumerate}

\end{enumerate}
\end{defi}

Returning to the general definition of a 2-groupoid, 
we emphasize that the first tensor product stands for the
operation of gluing surfaces (which should be thought of as
cobordisms between one dimensional manifolds), while the second
is induced by the first at the level of morphisms. 
In fact all  operations one  can imagine at the topological
level have natural counterparts in the groupoid setting. 
For instance capping off boundary circles with 
disks, or identifying  two boundary circles induce maps at the
homomorphism level. 
These 
maps correspond to the connected sum either with a disk or  with 
a cylinder, and thus   come from the tensor structure.

Another versions for the duality groupoid can be constructed by using 
only some of the possible gluings along boundaries, for instance 
by asking the common 
boundary contain only one circle, or
$\partial_-\Sigma=\partial_+\Sigma'$. 
In all these cases the presentation theorem below
has immediate reformulations, without introducing other generators or 
relations. 

Remark that one has  an embedding of the tower of mapping class groups
${\cal M}_{*,*}$, of  surfaces with boundary, in the
groupoid ${\cal D}$. This map associates to an element $\varphi$ of the 
${\cal M}(\Sigma)$ the element $(1,\varphi_*)
\in Hom((\Sigma,r),(\Sigma,r'))$,
where $\varphi_*$ transforms the rigid structure $r$ into
$ \varphi (r)$.

Observe also that all morphisms of the
groupoid  have representatives of the form $(1,c)$, and also that 
not all of these representatives come from elements of the mapping class group. 
In fact a necessary and sufficient condition for $(1,c)$ to be 
in the image of the mapping class group
is that the transformation  $c$ preserves the combinatorial 
configuration (i.e. the dual graph of the pants decomposition) 
of the rigid structure. In that case the rigid structures 
$r$ and $r'$ define uniquely (up to 
isotopy) a homeomorphism $\varphi$  such that $c=\varphi_*$.    

Let us explain what  the presentation of a
2-groupoid should be. Assume for simplicity that ${\cal O}(C)$ is an Abelian
category having  direct sums. 
\newtheorem{defi3}[defi]{Definition}
\begin{defi3}
A system of generators for the 2-groupoid $C$ consists of a 
collection of elements $x_i\in Hom(U_i,V_i)$, $U_i,V_i\in {\cal O}(C)$, $i\in
I$, such that:
\begin{enumerate}
\item Each $U\in {\cal O}(C)$ can be written as 
\begin{eqnarray*}
U=\bigoplus_{j=1}^n \bigotimes_{k=1}^{m_j} U_{i_{j k}}, \mbox{where}\:
i_{jk}\in I,\: n,m_j\in {\bf Z}.
\end{eqnarray*}
% $J\subset I^{{\bf Z}\times {\bf Z}}$.
\item Each $x\in Hom(U,V)$  can be written as   
\begin{eqnarray*}  
x= \oplus_{j} \otimes_{k} \circ _{l=0}^mx_{i_{jkl}}
\end{eqnarray*}
 where 
$\circ$ is the usual composition of morphism (subject to the 
source=target condition)
$U=\bigoplus_{j} \bigotimes_{k} U_{i_{jk0}}$,
$V=\bigoplus_{j} \bigotimes_{k} V_{i_{jkm}}$, and each of the 
$x_{i_{jkl}}$ is either equal to the identity morphism, 
or is one of the generators. 
\end{enumerate}
\end{defi3}

\newtheorem{defi4}[defi]{Definition}
\begin{defi4}
A presentation of a 2-groupoid $C$ is given by the system 
of generators $x_i\in Hom(U_i,V_i)$, and a system of
relations $r_j\in Hom(Z_j,W_j)$ that can be written 
in terms of the $x_i$'s.
\end{defi4}

The 2-groupoid with presentation $<x_i, i \mid r_j, j >$
can be constructed abstractly in the following way.
Fix the set of objects ${\cal O}(C)$.
 For  $U,V\in {\cal O}(C)$
define 
\[ Hom_0(U,V)= \bigoplus_{j}\bigotimes_{k}\circ_{l}
Hom_{00}(U_{i_{jkl}},V_{i_{jkl}}), \]
where $Hom_{00}(U_i,V_i)$ is the set of those maps constructed 
from  the $x_j$ with the same source and target.

Set $ RHom(U,V)\subset  Hom_0(U,V)$ for the subset 
of those $\varphi$ which can be written as $\alpha \circ \psi \circ
\beta $ with $\psi $ of the form
\[ \psi = \bigoplus_{j} \bigotimes_{k} \psi_{i_{jk}}, \]
where, for each $j$,  some of the  elements $\psi_{i_{jk}}$ 
are relations $r_l$ and the others are identity morphisms.
The set 
$Hom(U,V)=  Hom_0(U,V)/ RHom(U,V)$ is by definition the 
set of morphisms between $U$ and $V$. 

 The main purpose of this paper is to prove that the Moore-Seiberg
equations give a presentation of the 2-groupoid ${\cal D}$. 

\subsection{The geometric point of view}
Let us discuss an analogous situation. One can define
 a group presentation
$G=<x_i, i \mid r_j, j>$ geometrically as follows. Fix a
basepoint, and for each
generator $x_i$ a loop, then  attach a 2-cell on each loop made up from a word 
$r_j$. The  space $X_G=\bigvee_{x_i} S^1 \bigcup _{r_j} D^2$ has
the fundamental group equal to  $G$. 

Let us go one step further, to the presentation of a 1-groupoid $C$. 
 Consider a presentation of 
a 1-groupoid $C$ given by $<x_i, i \mid r_j, j>$, where 
$x_i\in Hom(s(x_i),t(x_i))$,  $s, t$ being  the source and target 
maps. Construct  a 2-complex $X_C$ in the same vein, by identifying
the set $F$ of final objects with a set of $0$-cells
and by choosing a $1$-cell  
connecting  $a$ and  $b$ in $F$, for each $x_i$ 
such that  $s(x_i)=a$, $t(x_i)=b$. Attach  a 2-cell on 
a loop representing $r_j$, for each $j$. The fundamental  groupoid 
$\pi_1(X_C, F)$ with base points in $F$ is the 1-groupoid of the given
presentation. 
Notice that relation $r_j$ with $s(r_j)\neq t(r_j)$  add further
identifications in $F$,  to enable us to attach 
 2-cells. 

Consider now a 2-groupoid $C$, with 
generators and relations $x_i$ and $r_i$. 
Like before, identify the final objects
of ${\cal O}(C)$ with the set of
$0$-cells, and add a $1$-cell between $s(x_i)$  and $t(x_i)$ 
for each generator  $x_i$. 
Next let   $K^1$ be the 1-complex obtained as closure of this 
structure under the tensor product, 
meaning that each edge $x_i$ induces the attachment of 
other edges, denoted by $x_i\otimes 1_a$ (respectively 
 $1_a\otimes x_i$), between 
$s(x_i)\otimes a$ and $t(x_i)\otimes a$. These  correspond to  
elements $x_i\otimes 1_a\in Hom(s(x_i)\otimes a,t(x_i)\otimes a)$. 
Recall that $1_a$ is the identity element in the 
group $Hom(a,a)$. 

Attach  to $K^1$ 2-cells along  the loops associated to the relations 
$r_j$, and take the $\otimes$-closure $K^2$, meaning that once a 
2-cell is attached on the vertices $u_i$ and edges $e_i$ then all its 
translated copies on the vertices $u_i\otimes a$ and edges $e_i\otimes
1_a$ (respectively $a\otimes u_i $ and $1_a\times e_i$) are also
2-cells.  Finally, add the DC-cells that 
come from the tensor structure. These cells are defined as follows. Assume 
that we have $a\in Hom(x, x')$ and $y\in Hom(y,y')$. Consider the 
four vertices $x\otimes y$, $x'\otimes y$, $x\otimes y'$ and 
the edges $a\otimes 1_y$, $1_{x'}\otimes b$, 
$1_x\otimes b$ and $a\otimes 1_{y'}$ relating these vertices. 
Attach a 2-cell on the 
square made off the edges and call it a DC-cell (disjoint
commutativity). The relations expressed by these cells are the 
obvious $(a\otimes 1_y)(1_{x'}\otimes b)= (1_x\otimes b)(a\otimes
1_{y'})$. Call $X_C$ the new $2$-complex. 
The tensor multiplication gives a multiplicative structure on 
the groupoid of paths in $X_C$. 
When adding the 2-cells one obtain a tensor multiplication on the 
fundamental groupoid $\pi_1(X_C,{\cal O}(C))$, and the
2-groupoid obtained this way is isomorphic to $C$.

%The 
%Harer-Hatcher-Thurston technique of proving that a certain group 
%arising from a topological construction has a given  presentation
%can be adapted to $2$-groupoids in the following way.
%Consider the 
%2-complex $X_P$ associated to a given presentation $P$ of $C$. 
%after we
%verified that the algebraic relations hold indeed in the groupoid $C$. 
%We add more 2-cells corresponding to all relations which we know to 
%hold in $C$: the generators are
%no more free generators but considered as elements of $C$. We obtain
%then the 2-complex $X_{P,C}$. If we want to identify $C$ with the 
%2-groupoid with the presentation $P$ it suffices then to prove that 
%$X_{P,C}$ is connected and simply connected. In order to see that 
%there are not redundant relations, we have to show that once a
%relation is deleted from $P$, the associated 2-complex is not simply 
%connected.  

\section{The Moore-Seiberg equations}
\subsection{Main results}

This section contains the main results of the paper.
 
\newtheorem{theo}{Theorem}[section]
\begin{theo}
The duality  2-groupoid  ${\cal D}$ has the  2-groupoid presentation with: 
\begin{enumerate}
\item[] Generators $T_1, R, B_{23}, F, S, A, D, P$ and their inverses. 
\item[]  Relations (Moore-Seiberg equations)\\
1. at the level of a pair of pants:\\
\ a) $T_1B_{23}=B_{23}T_1$,\ $T_2B_{23}=B_{23}T_3$, \ $T_3B_{23}=B_{23}T_2$,
where $T_2=RT_1R^{-1}$ and $T_3=R^{-1}T_1R$,\\
\ b) $B_{23}^2=T_1T_2^{-1}T_3^{-1}$,\\
\ c) $R^3=1$,\\
\ d) $RB_{23}R^2B_{23}RB_{23}R^2=B_{23}RB_{23}R^2B_{23}$,\\
2. relations defining  inverses:\\
\ a) $P^{(12)}F^2=1$,\\
\ b) $T_3^{-1}B_{23}^{-1}S^2=1$,\\
3. relations coming from ``triangle singularities'':\\
\ a) $P^{(13)}R^{(2)}F^{(12)}R^{(2)}F^{(23)}R^{(2)}
F^{(12)}R^{(2)}F^{(23)}R^{(2)}F^{(12)}=1$,\\
\ b) $T_3^{(1)}FB_{23}^{(1)}FB_{23}^{(1)}FB_{23}^{(1)}=1$,\\
\ c) $B_{23}^{-1}T_3^{-2}ST_3^{-1}ST_3^{-1}S=1$,\\
\ d) $R^{(1)}(R^{(2)})^{-1}FS^{(1)}FB_{23}^{(2)}
B_{23}^{(1)}=FS^{(2)} T_3^{(2)}(T_1^{(2)})^{-1}B_{23}^{(2)}
F$,\\
4. relations coming from the symmetric groups:\\
\ a) $P^2=1$, \\
\ b) $P^{(23)}P^{(12)}P^{(23)}=P^{(12)}P^{(23)}P^{(12)}$.\\
5. relations involving annuli and disks:\\
\ a) $ A^{(12)} D^{(23)} = A^{(23)} D^{(12)}$, \\
\ b) $A^{(12)}D_3^{(13)}=
A^{(13)}D_3^{(13)} F$,\\
\ c) $A^{(12)}A^{(23)}=
A^{(23)}A^{(12)}$ \\
\ d) $S D = D S$. \\
\end{enumerate}
\end{theo}

We used the convention that superscripts tell us
on which factors of the tensor product the move acts. Here the  tensor
structure is implicit.

It was  proved in 
\cite{Fun} that any topological invariant of 3-manifolds determines a 
unique  maximal associated
TQFT. In the terminology of \cite{Wa} and \cite{Gel} this is a
 TQFT with corners. Notice that a  TQFT gives rise to a representation 
$\rho_{*}:{\cal M}_{*}\longrightarrow End(W_*)$.

When considered on a torus, that is when capping the
1-holed torus with a disk, relations 2.b) and 3.c) 
give rise to the well known  morphism from $SL(2,C)$ into the groupoid
of moves acting on the torus, which groupoid contains  
the mapping class group of the torus, as a maximal group.
If on a sphere with four holes we factor  out by the twists around the
holes, i.e. if  we consider the groupoid of moves on a 
fourth punctured sphere, then relations 2.a) and 3.b) give rise
to another morphism from $SL(2,C)$ into the groupoid of the
4-holed sphere. This latter morphism is used in the 
classification of $2$-bridge knots.

\newtheorem{theo3}[theo]{Theorem}
\begin{theo3} 
Assume that the TQFT is finite dimensional, unitary, cyclic and has a
unique vacuum (see section 5 for complete definitions).
Then  $\rho_*$ extends canonically to a representation of the 
the full duality groupoid ${\cal D}$. 
In particular all maximal TQFT representations verify the 
Moore-Seiberg equations. 
\end{theo3}
%Observe that, conversely a representation of ${\cal D}$ 
%verifying the Moore-Seiberg equations   is the extension 
%of the monodromy representation of a TQFT in dimension 3. 
%This is more commonly stated as the fact that a conformal field theory
%(CFT) in dimension 2 gives naturally a 3-dimensional TQFT. 
%This result is due to  Witten \cite{Wi}, explicitly stated by
%Kontsevich and reproved in many contexts by L.Crane, P.Degiovanni and
%others. 
On the other direction  Kohno (see \cite{Kohno,Kohno2}) used the data coming from 
conformal field theory to construct representations of the tower of
mapping class groups (and in fact of the duality groupoid). He proves then
that these determine
topological invariants for 3-manifolds (which actually extend to a TQFT).     

\subsection{Topological interpretation of the generators and
  relations}

The generators written above  can be explicitly viewed in
the topological picture of the groupoid, so that the relations 
become tautological. We have: 

\begin{enumerate}
\item Moves on  rigid structures on a pair of pants, which are the 
three twists $T_j$ around the boundary circles, the 
knob twisting  $B_{23}$ and the  cyclic  permutation $R$ of the numbering
of the boundary components. 
%\item A twist or a move of type $RC$ on a cylinder, which can be
%  obtained from $R$ by capping off the boundary component of the 
%trinion labeled 3.
\item The move $F$ on a sphere with 4 holes decomposed into
two pairs of pants. The decomposition curve is transformed
as such that the new curve does not intersect the seams that 
the old one did, and intersects each of the other two seams
exactly once. The numberings of the pairs of pants transform
such that the decomposition curve remains labeled by $1$ and the
boundary curve labeled by $2$ of the first pair of pants becomes the
the curve labeled by $3$ of the first pair of pants.  
\item The move $S$ on the 1-holed torus. $S$ changes the 
rigid structure as the element 
$\left [\begin{array}{cc}
0 & -1 \\
1 & 0 
\end{array}
\right ]$ of the mapping class group. 
\item The move $P$ which transposes  the numberings 
of two pairs of pants, two annuli or two disks. Using the groupoid
structure, $P$ generates the whole  permutation group of 
numberings.  
\item Moves $D$ and $A$ which correspond to contracting annuli and disks.
 Their inverses consist of expansions of disks or
annuli.
\end{enumerate}

\begin{figure}
\centering
\leavevmode
\epsfxsize=3in
\epsfysize=1.2in
\epsfbox{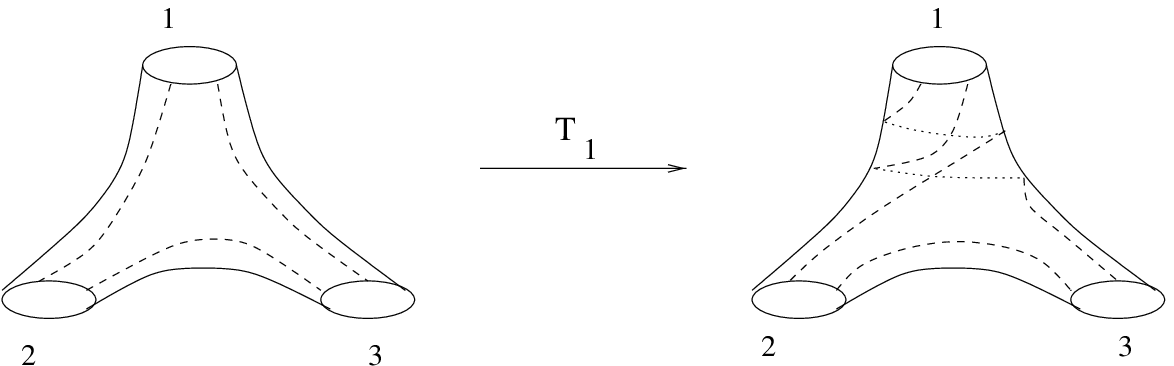}
\caption{$T_1$}
\label{T1}
\end{figure}

\begin{figure}
\centering
\leavevmode
\epsfxsize=3in
\epsfysize=1.2in
\epsfbox{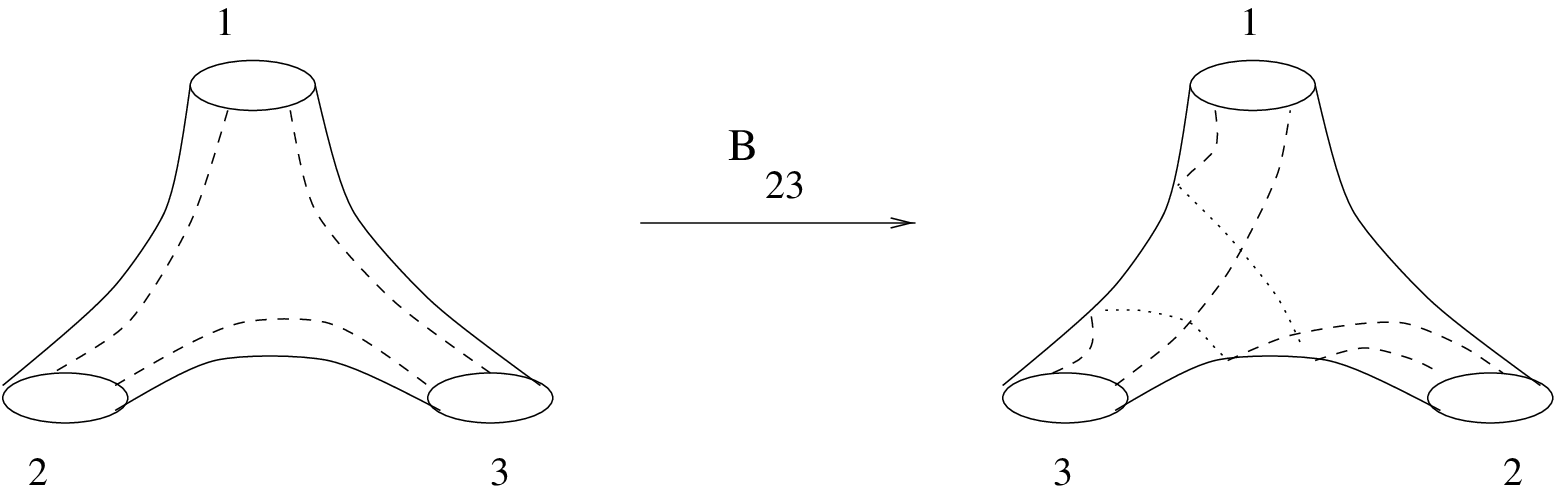}
\caption{$B_{23}$}
\label{B23}
\end{figure}

\begin{figure}
\centering
\leavevmode
\epsfxsize=3in
\epsfysize=1.2in
\epsfbox{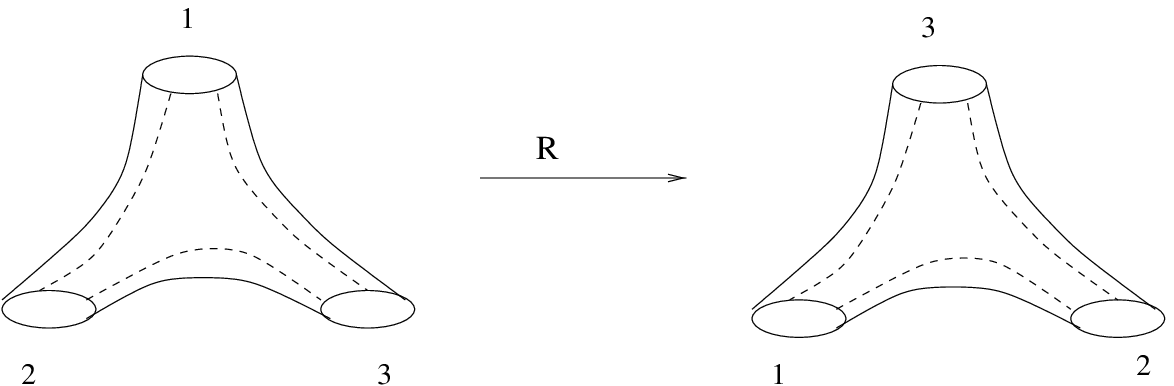}
\caption{$R$}
\label{R}
\end{figure}

\begin{figure}
\centering
\leavevmode
\epsfxsize=3in
\epsfysize=1.2in
\epsfbox{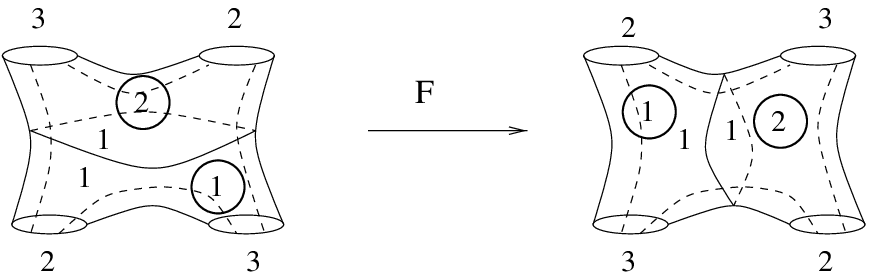}
\caption{$F$}
\label{F}
\end{figure}

\begin{figure}
\centering
\leavevmode
\epsfxsize=3in
\epsfysize=1.2in
\epsfbox{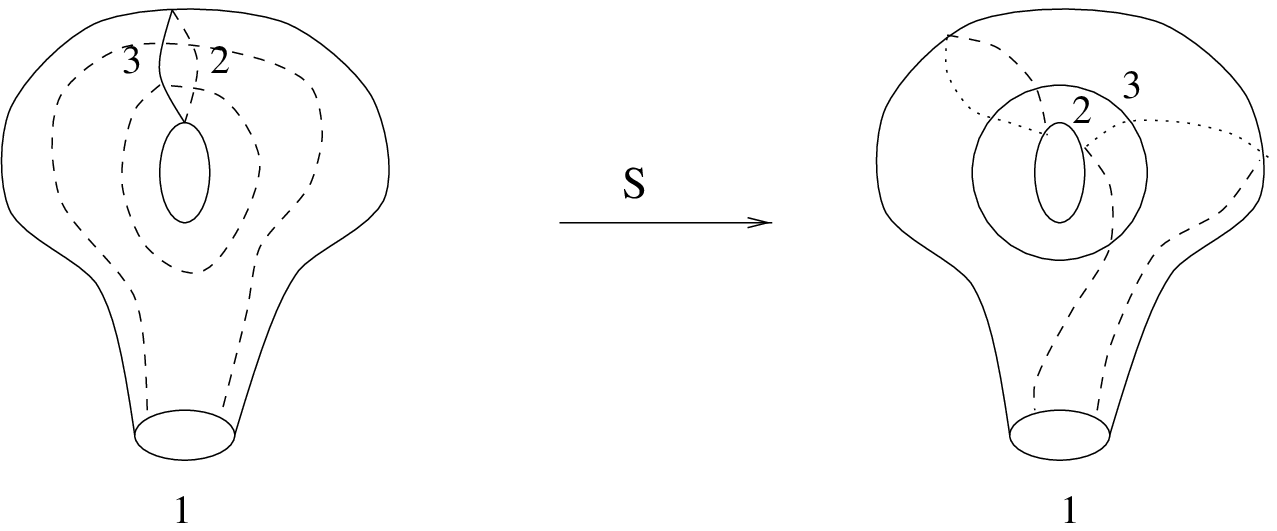}
\caption{$S$}
\label{S}
\end{figure}

These elementary moves are described in Fig.~\ref{T1}-\ref{S}. 
In these figures the convention is   that the circles 
of the DAP-decomposition are drawn as plain curves, every curve 
being labeled by a number in each elementary surface that it bounds,
the seams are pictured as dashed curves, and  each pair of pants  
carries an encircled number, these numbers defining the ordering of
elementary surfaces. If one
element of this data is absent this means that it can be chosen
 arbitrarily in the given situation. 
 Note  that the relations given in Theorem 3.1 can be easily 
verified pictorially.

Let us stress  out that moves represent changes of rigid 
structures and not homeomorphisms. 
 The first group of relations
  are identical with the ones giving the 
presentation of ${{\cal M}}_{0,[3]}$, (the extended mapping class group of the 
3-holed sphere, in which is allowed to interchange 
the boundary components).
However, the moves $F$, $A$ and $D$ do not have analogues at the
level of homeomorphisms.
%The relations defining inverses are natural too: if we forget the 
%markings then $F^2$ is simply the identity map. However once we add
%the labelings of the plain curves then the map $F^2$ acts like 
%a permutation on the numberings. On the other hand the map $S^2$ 
%can be expressed in terms of the generators of the trinion, the 
%holed torus being viewed as a pair of pants, two boundary circles of 
%which being identified. 

%\begin{figure}
%\centering
%\leavevmode
%\epsfxsize=3in
%\epsfysize=1.2in
%\epsfbox{Newfig/B12.eps}
%\caption{$B_{12}$}
%\label{B12}
%\end{figure}

%\begin{figure}
%\centering
%\leavevmode
%\epsfxsize=3in
%\epsfysize=0.3in
%\epsfbox{Newfig/RC.eps}
%\caption{$RC$}
%\label{RC}
%\end{figure}

\subsection{The 2-complex $\Gamma$}

For the proof of Theorem 3.1 we will adapt the Harer-Hatcher-Thurston technique
to the present situation.
 To this end we construct a family of  2-complexes,  related by
the tensor  product. 

\newtheorem{ull}[theo]{Definition}
\begin{ull}
The complex $\Gamma(\Sigma)$ is obtained as follows:
\begin{enumerate}
\item Its vertices are the various rigid structures on $\Sigma$.
\item Two vertices are related by an edge if there is one
  transformation of type $B_{23}, T_1, R, F, S, P, D, A$ 
which relates the respective rigid structures. 
\item The first set of 2-cells  is given by the Moore-Seiberg 
equations: each equation gives a circuit on the 1-skeleton 
and we attach a 2-disk on it. 
\item The second set of 2-cells are the DC-cells which 
represent the commutation relation between two moves whose 
supports are sub-surfaces with disjoint interiors. 
\item The third set of relations correspond to relations among the
permutations in the ordering of the elementary surfaces in
the DAP-decomposition.
\end{enumerate}
\end{ull}

We observe that when talking about the moves $F, S,...$ we already make use 
 of the tensor structure on the duality groupoid, because these
moves are defined  on sub-surfaces and are extended by identity
outside the support. Hence it makes sense  to consider the 
union $\Gamma= \bigcup_{\Sigma} \Gamma(\Sigma)$. The set of vertices
 has a multiplicative  structure, the tensor product of  the groupoid, and 
the fundamental groupoid $\pi_1(\Gamma)$ is nothing but the 
2-groupoid with presentation given by  the Moore-Seiberg equations. 
Thus Theorem 3.1 follows from   

\newtheorem{uli}[theo]{Theorem}
\begin{uli}
The complex $\Gamma(\Sigma)$ is connected and simply connected. 
The mapping class group ${\cal M}(\Sigma)$ acts freely on it. 
\end{uli}

The proof is reminiscent of \cite{HaTh}. We consider first simpler 
structures which mimic the construction of $\Gamma(\Sigma)$.
Thus  we start
with the groupoid of markings, then add overmarkings and eventually 
come to the last object. The cases of markings
and  overmarkings are  solved with techniques of Cerf
theory, and  simple algebraic topology arguments yield the result for
rigid structures. The proof of this result will be done in detail in 
Chapter 4.

\section{Proof that Moore-Seiberg equations are sufficient}
\subsection{Elements of Cerf theory}

Following \cite{HaTh}, given a surface
 we  call {\em marking} a finite collection
of disjoint simple closed curves lying in the interior of the surface,
 that decompose the  surface into pairs
of pants. Thus markings are obtained from rigid structures by forgetting
 the  annuli, disks,  seams and numberings.
Of course for a surface to admit a marking it must be different
from a sphere, torus, disk or annulus. 

Let $\Sigma$ be a given surface. We want to find  a presentation
of the groupoid of transformations of markings on  $\Sigma$.
For this we will use Cerf theory \cite{Cerf},  in
an analogous way  it was used in \cite{HaTh} for the study of
the mapping class group.

For each marking there exists a Morse height function $f:\Sigma\rightarrow
{\bf R}$, such that the decomposition curves are
connected components of level sets of $f$. 
The space ${\cal F}$ of height functions has a stratification
\begin{eqnarray*}
{\cal F}={\cal F}^0\cup {\cal F}^1\cup {\cal F}^2 \cup {\cal F}^3\cup \cdots
\end{eqnarray*}
where ${\cal F}^k$ are strata of codimension $k$.

In particular ${\cal F}^0$ is the set 
  of Morse functions
(i.e. functions with finitely many critical points, all Morse
and at different heights), and is a dense open subset of ${\cal F}$,
and $ {\cal F}^1={\cal F}^1_{\alpha}\cup{\cal F}^1_{\beta}$, with
$ {\cal F}^1_{\alpha}$ having finitely many critical points, all Morse
and at different heights, except for   one which is a birth-death 
(i.e. the height function looks locally like $f(x,y)=\pm x^3\pm y^2$),
and ${\cal F}^1_{\beta}$ having finitely many critical points, all Morse
and at different heights, except for a pair of Morse points which 
are at the same height.  

We will make use of the two results results in Cerf theory given in the sequel.

\newtheorem{theo4}[theo]{Theorem}
\begin{theo4} 
 Any two Morse functions $f_0$ and $f_1$
can be joined by a path of height functions $(f_t)_{t\in[0,1]}$,
with the property that all $f_t$ are Morse except for finitely
many, and for these exceptional functions
 the path crosses ${\cal F}^1$ transversely.
\end{theo4}

\begin{figure}[htbp]
\centering
\leavevmode
\epsfxsize=2in
\epsfysize=1in
\epsfbox{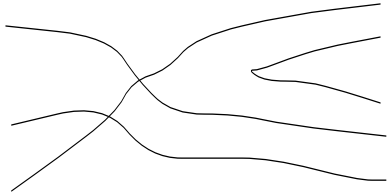}
\caption{Example of a graph}
\label{graf}
\end{figure}

A path having the property described in the theorem 
is called a good path.
For a better understanding it is customary to sketch the graph of
a path, namely to trace the critical points of the functions 
$f_t$, $t\in [0,1]$. An example of a  graph for  a good path 
 is given in Fig.~\ref{graf}.

The second theorem tells us how a homotopy of paths crosses the
codimension two stratum.

\newtheorem{theo5}[theo]{Theorem}
\begin{theo5} 
If $(f_t)_t$ is a closed good path in ${\cal F}$, then
there exists a homotopy $(f_{t,u})_{u}$  from it to the constant
path such that $f_{t,u_0}$ is good for every $u_0$, except for the
following isolated singularities:
\begin{itemize}
 \item[a).] crossing  (two crossings are cancelled
or introduced),
\item[b).] birth-death  (two birth-death points are canceled or 
introduced), 
\item[c).] triangle  (three non-degenerate critical points 
lie at the same level),
\item[d).] beak singularity  (a birth-death point crosses a 
non-degenerate critical 
point),
\item[e).] swallow tail,  
\item[f).] two birth-death or crossing points occur simultaneously.
\end{itemize}
 \end{theo5}

\begin{figure}[htbp]
\centering
\leavevmode
\epsfxsize=3in
\epsfysize=3.2in
\epsfbox{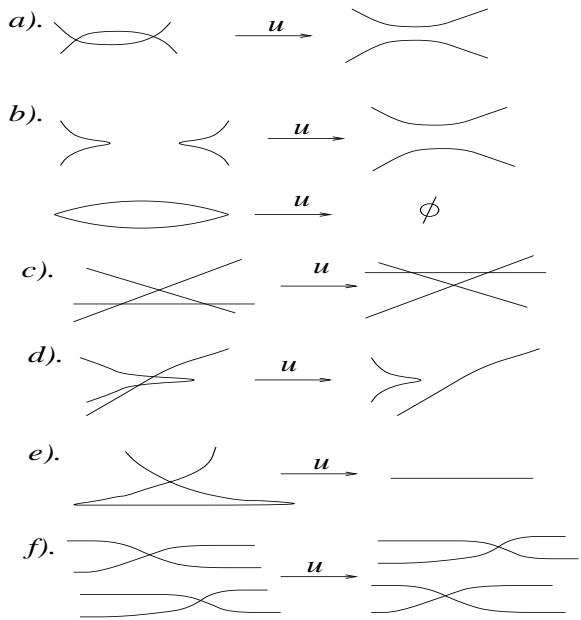}
\caption{Singularities}
\label{sing1}
\end{figure}

These singularities are shown graphically in Fig.~\ref{sing1}.

Given two markings, associate to them Morse functions that
determine them. One can pass from one function to the other along a good path.
 The only case when the marking can change is when one crosses the
codimension one stratum. The marking does not change at a beak point, nor
at a crossing point if any of the critical points that cross
has index different from  $1$. Hence the 
only interesting points are the crossings of saddle  points.
To understand what happens in this case, let us restrict
our attention  to the semi-local picture containing these points,
namely to the pairs of pants determined by the marking that
contain the two points. 
If the pairs of paints are disjoint, then the marking remains unchanged
after the crossing. 

If the pairs where the crossing occurs share some boundary
components, we have the situations described in Fig. 8.
Here and below to encode the crossings, rather than using the 
associated trivalent graph, as it was done in \cite{HaTh}, we will
use the ascending-descending  manifold model, which
is more suggestive in this situation. Let us recall that
the ascending (unstable) manifold is the submanifold on which
the quadratic function that gives the local model of
the singularity is positive definite, and the descending (stable)
manifold is the submanifold on which the quadratic function
is negative definite. In the case of index one singularities
on a surface both these submanifolds are one-dimensional.
On the left hand column of Fig. 8 we represented the 
descending  manifold model viewed from above, and
one should imagine the two descending manifolds
exchanging heights when the crossing occurs.

\begin{figure}[htbp]
\centering
\leavevmode
\epsfxsize=4in
\epsfysize=4in
\epsfbox{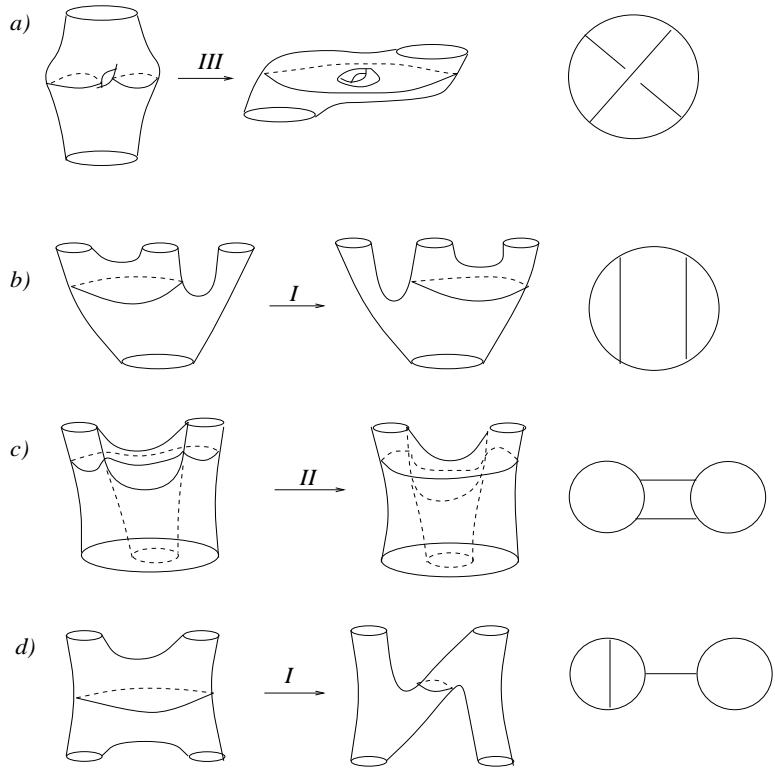}
\caption{Crossings}
\label{cross}
\end{figure}

These four types of crossings give rise to the four moves
of Hatcher and Thurston \cite{HaTh}. Recall that $IV$ is obtained by
capping off one boundary component of the torus by using a disk, in
the move $III$.
Consequently any two markings can be transformed one into the other by
 applying finitely many moves like these.

To find the relations that these moves satisfy, we will rely on the second 
theorem.
Since beak points do not interfere with markings, 
the only singularities that produce relations between moves
are a), which expresses the fact that each move
is its own inverse, f), which gives the disjoint commutativity
between moves that occur far away from each other,
and the triangle singularity.
The latter produces the most interesting relations, 
which we will describe below.

\subsection{The groupoid of markings}

\begin{figure}[htbp]
\centering
\leavevmode
\epsfxsize=5.8in
\epsfysize=4in
\epsfbox{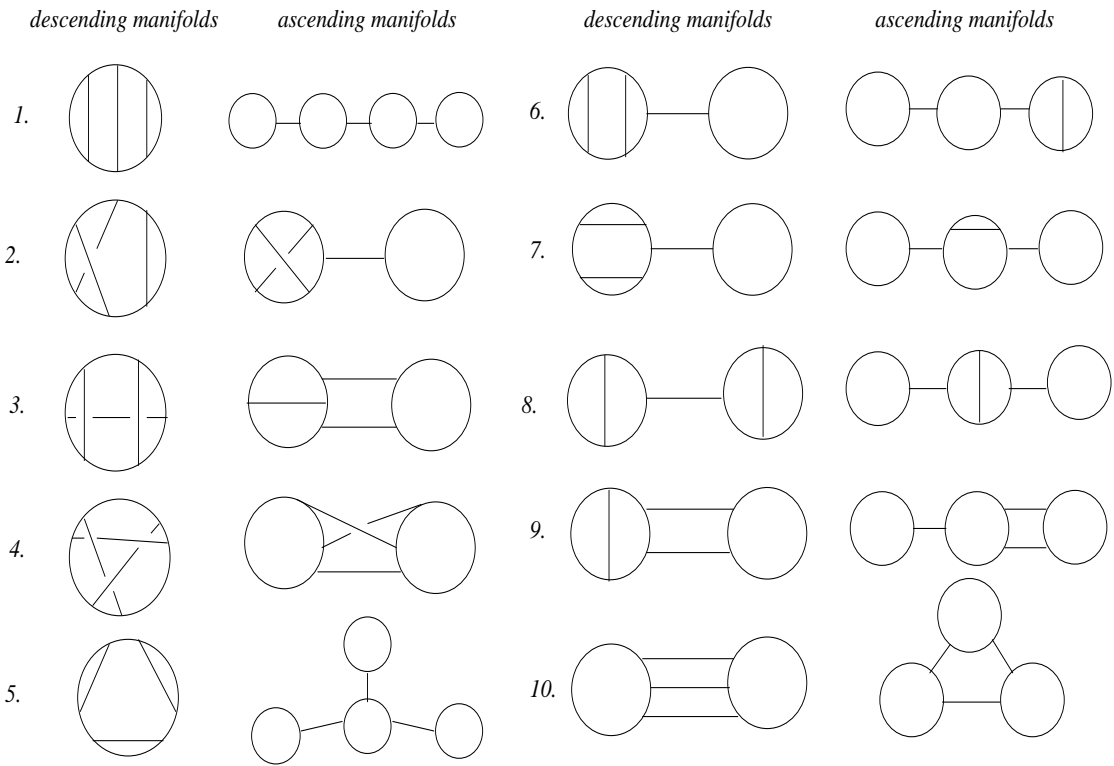}
\caption{Triangle singularities}
\label{config}
\end{figure}

Looking at the combinatorics of the circles below the singularity, 
and of the arcs connecting them, determined by the descending manifolds,
there are $20$ possible configurations. These configurations are described
in Fig.~\ref{config}. In   this figure  the second column consists of
 descending  manifold models which by changing $f$ to $-f$ are 
the ascending manifold models corresponding to the descending
manifold models from the first column.
Because of this symmetry, there are only $10$ relations 
between Hatcher-Thurston
moves arising from these singularities.
In these pictures the descending  manifolds are at 
different heights and the exchanges in heights correspond 
to crossings. A particular choice of heights specifies the 
vertex at which one begins tracing the boundary of the cell,
thus it suffices to make one choice for each diagram.

\begin{figure}[htbp]
\centering
\leavevmode
\epsfxsize=3.5in
\epsfysize=3.5in
\epsfbox{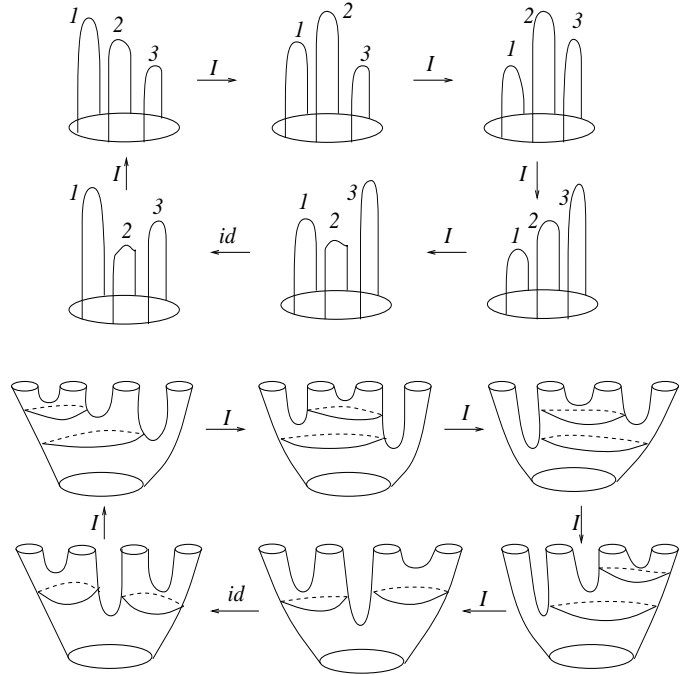}
\caption{Relation 1}
\label{cell1}
\end{figure}

It is not hard to see that
the configurations $1)$, $6)$, $7)$ and $8)$ give rise to the same relation.
So there are seven distinct relations 
coming from triangle singularities.
They are described  in Fig.~\ref{cell1}--~\ref{cell10}. For
clarity let us stress out
that relations 1,5,9 and 10 hold on a sphere with five holes,
while relations 2,3 and 4 hold on a torus with three holes. 

In the CW-complex setting, 
let us  consider  the 
 2-complex $\tilde{\Gamma}_0(\Sigma)$ defined as follows. 
  The vertices of  $\tilde{\Gamma}_0(\Sigma)$ are all 
 possible  markings on the surface  $\Sigma$,  and  
there is an edge 
between two vertices if there exists a transformation of type
$I, II, III$ or $IV$  relating the respective markings.  
The first set of 2-cells are the seven types of cells described 
in Fig.\ref{cell1}--~\ref{cell10}. To these we add the
cells that express disjoint commutativity, called DC-cells, which come
from crossing singularities.

\newtheorem{proposition}{Proposition}[section]
\begin{proposition}
The 2-complex $\tilde{\Gamma}_0(\Sigma)$ is connected and 
simply connected. 
\end{proposition}
{\em Proof:} This 
is a consequence of the two theorems from Cerf theory we 
cited previously, and the geometric interpretation given to 
markings.$\Box$  

\begin{figure}[htbp]
\centering
\leavevmode
\epsfxsize=3.5in
\epsfysize=3.5in
\epsfbox{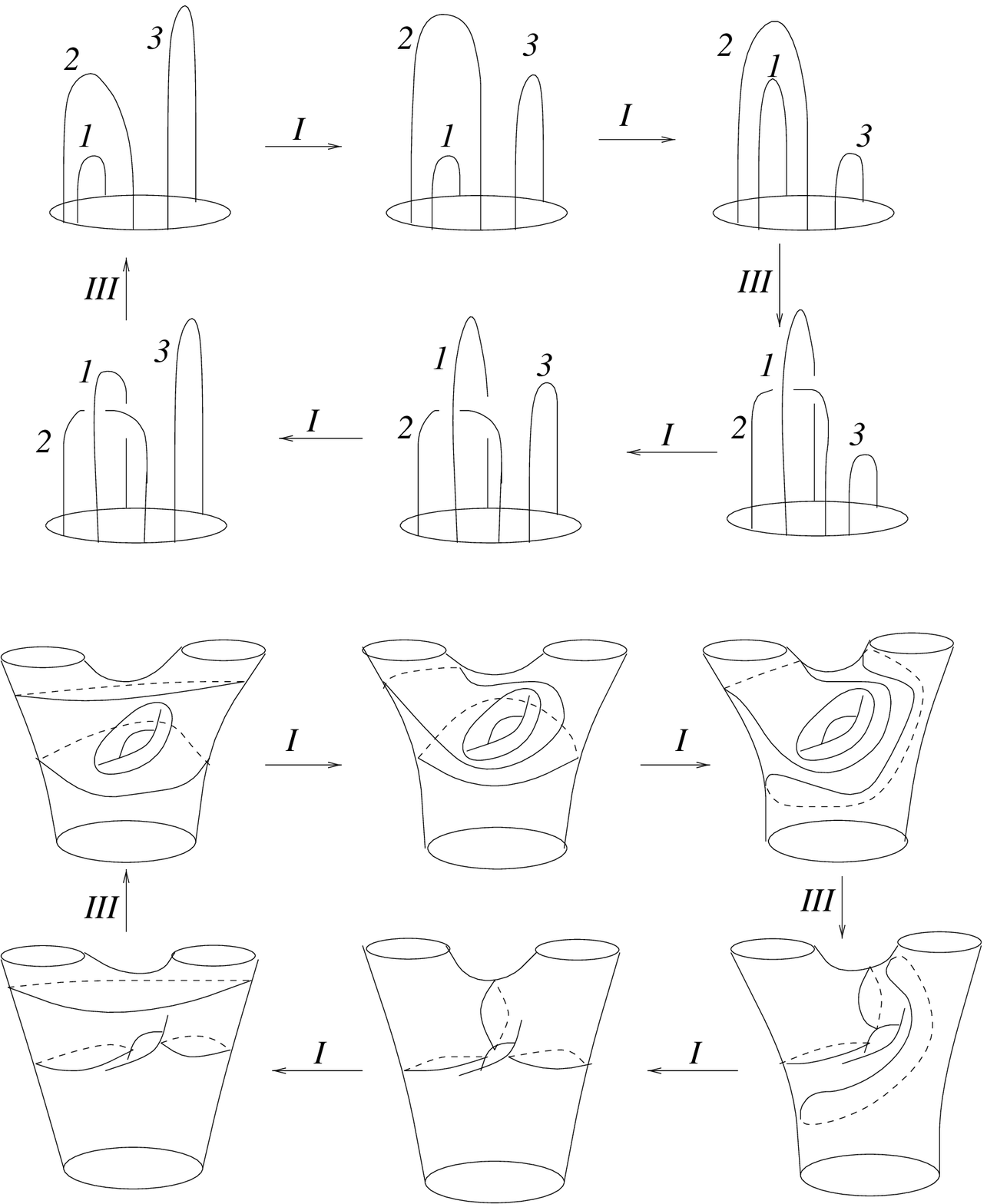}
\caption{Relation 2}

\centering
\leavevmode
\epsfxsize=4in
\epsfysize=4in
\epsfbox{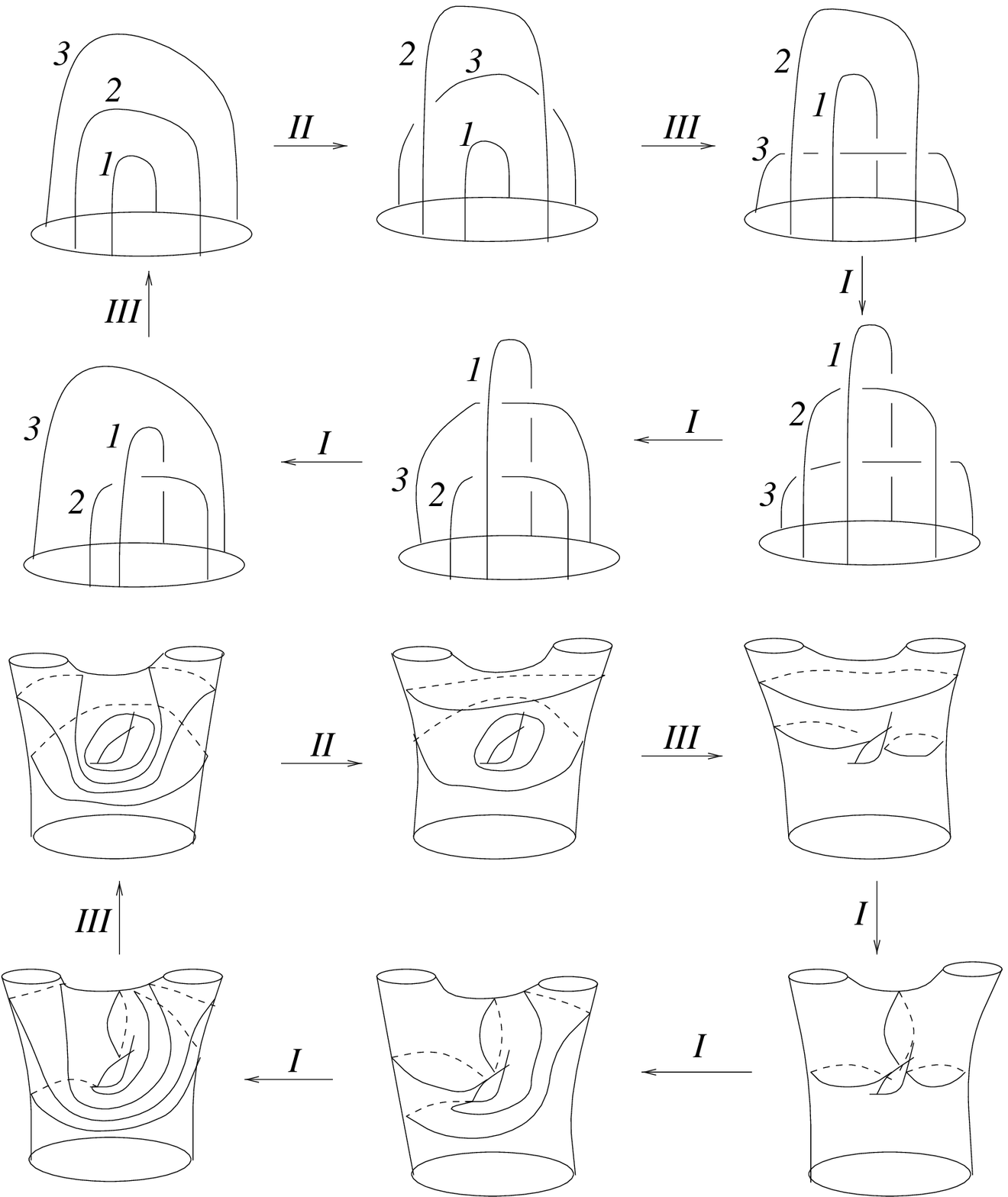}
\caption{Relation 3}
\end{figure}

\begin{figure}[htbp]
\centering
\leavevmode
\epsfxsize=4in
\epsfysize=4in
\epsfbox{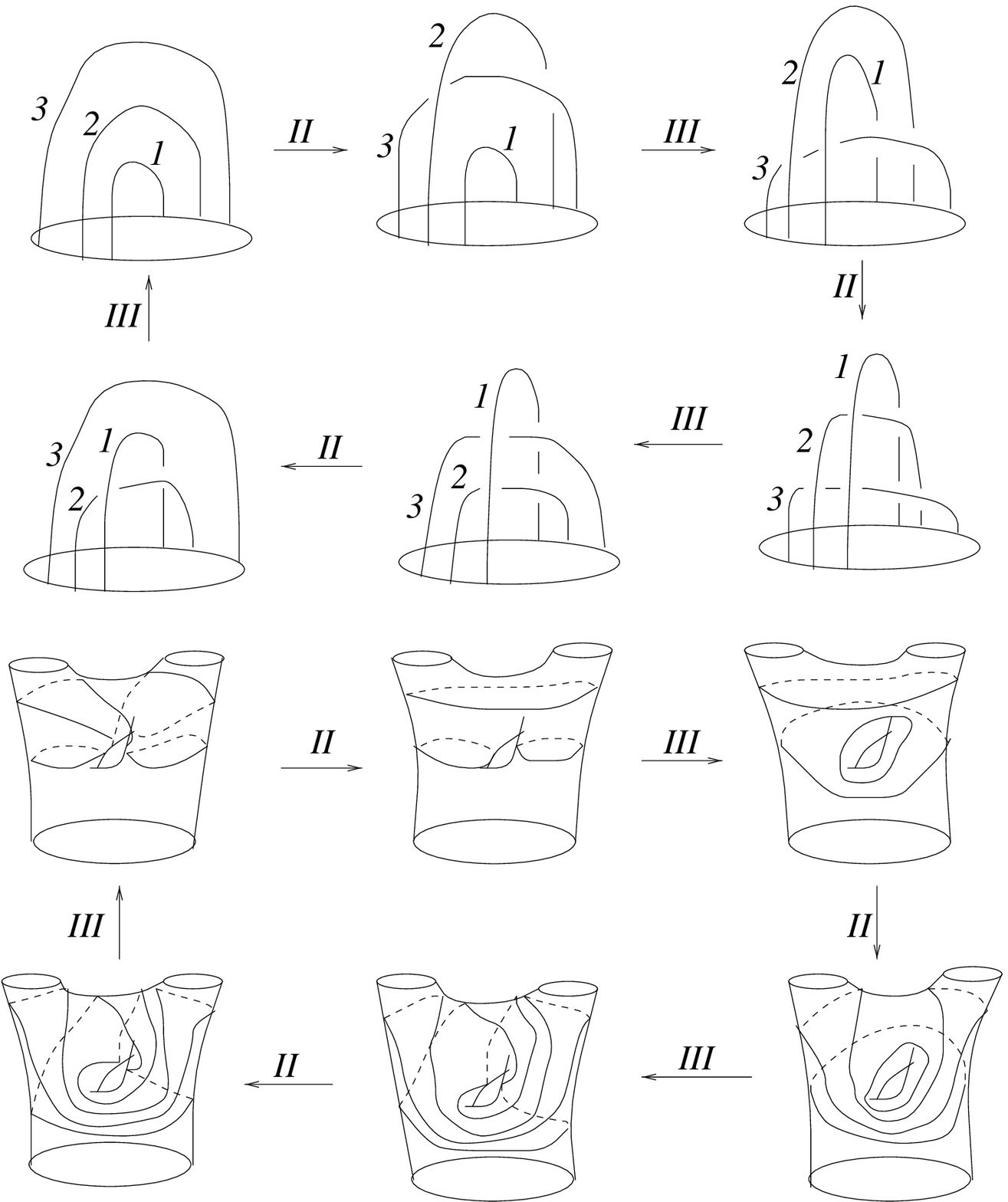}
\caption{Relation 4}

\centering
\leavevmode
\epsfxsize=4in
\epsfysize=4in
\epsfbox{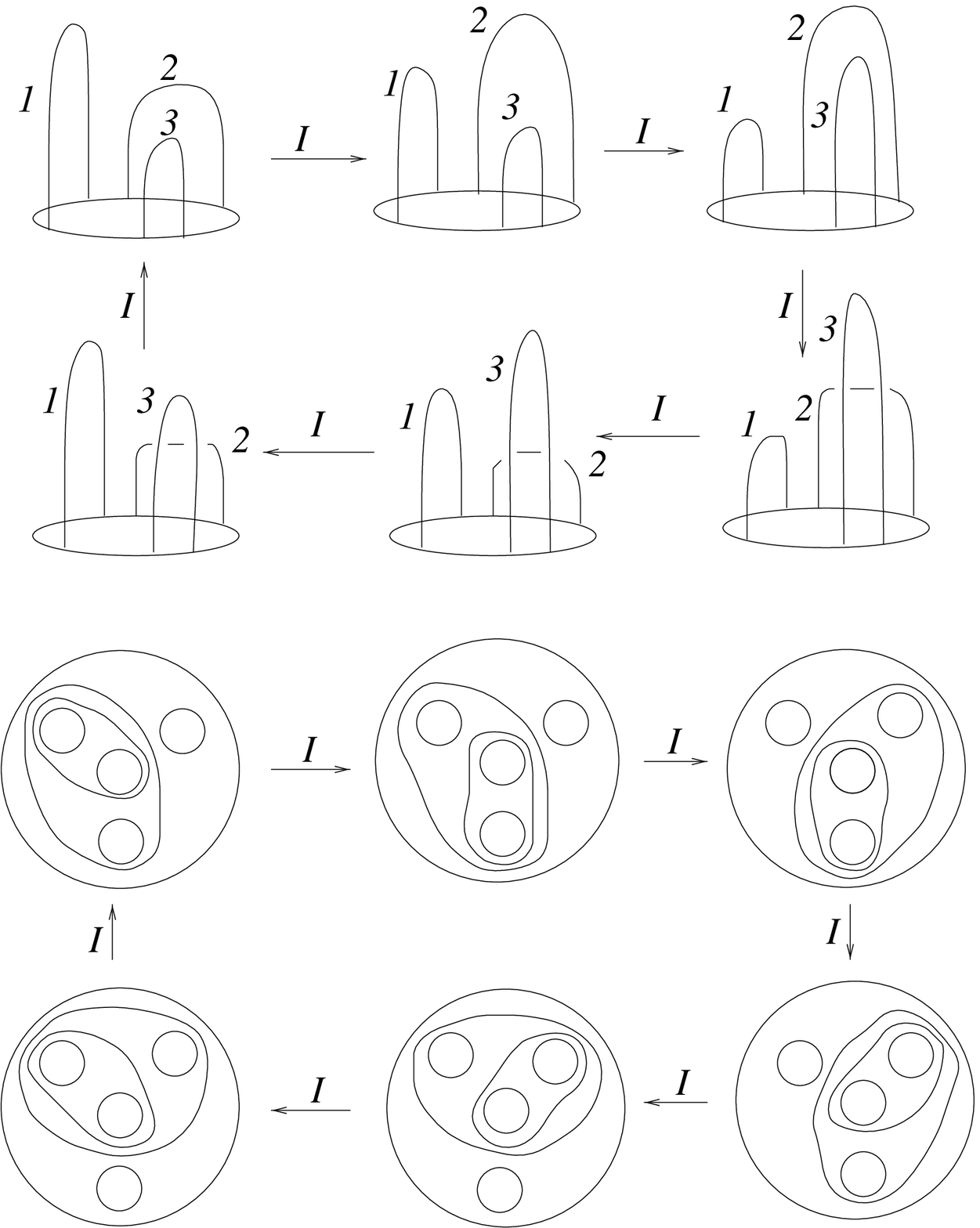}
\caption{Relation 5}
\end{figure}

\begin{figure}[htbp]
\centering
\leavevmode
\epsfxsize=4in
\epsfysize=4in
\epsfbox{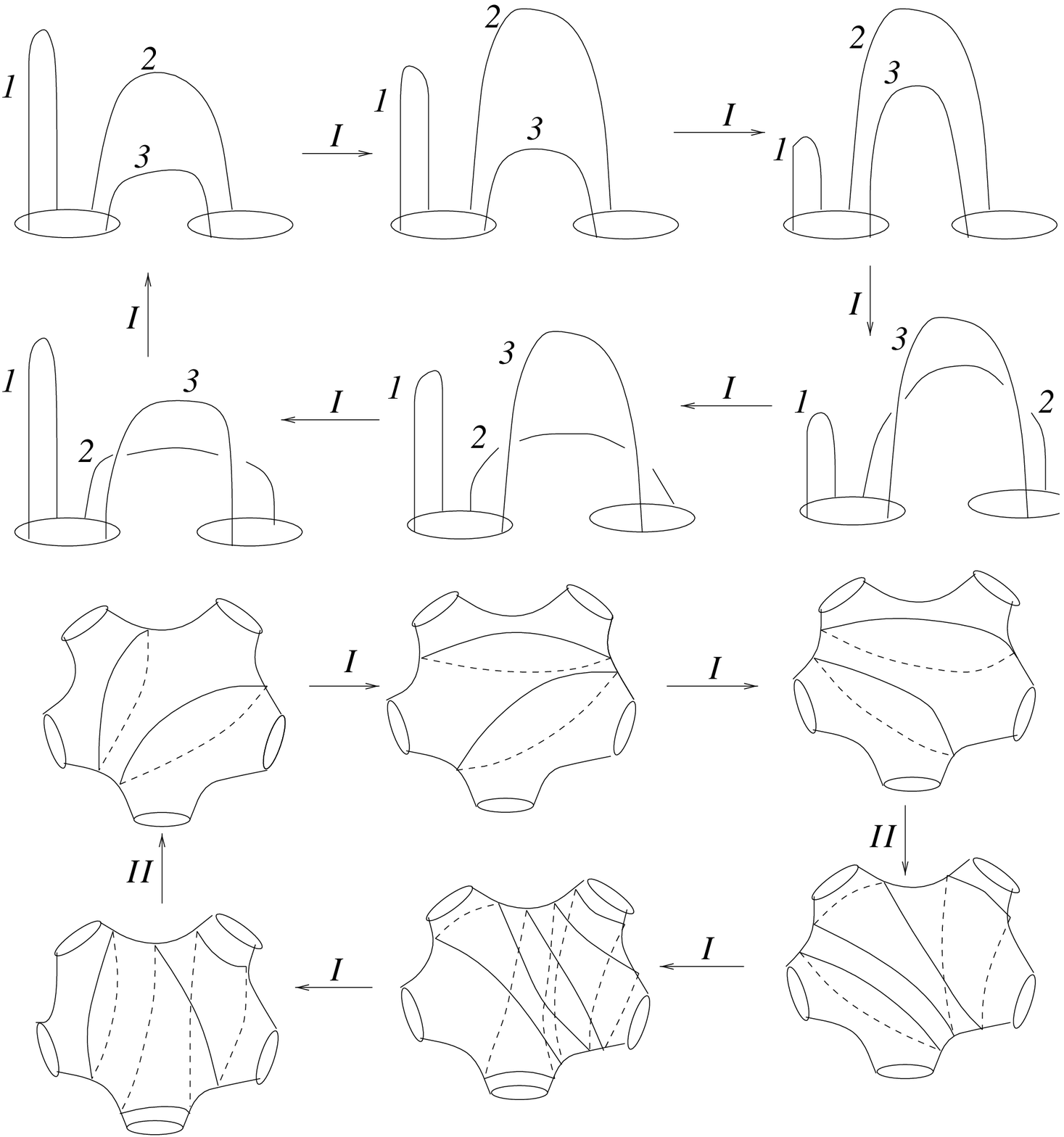}
\caption{Relation 9}

\centering
\leavevmode
\epsfxsize=4in
\epsfysize=4in
\epsfbox{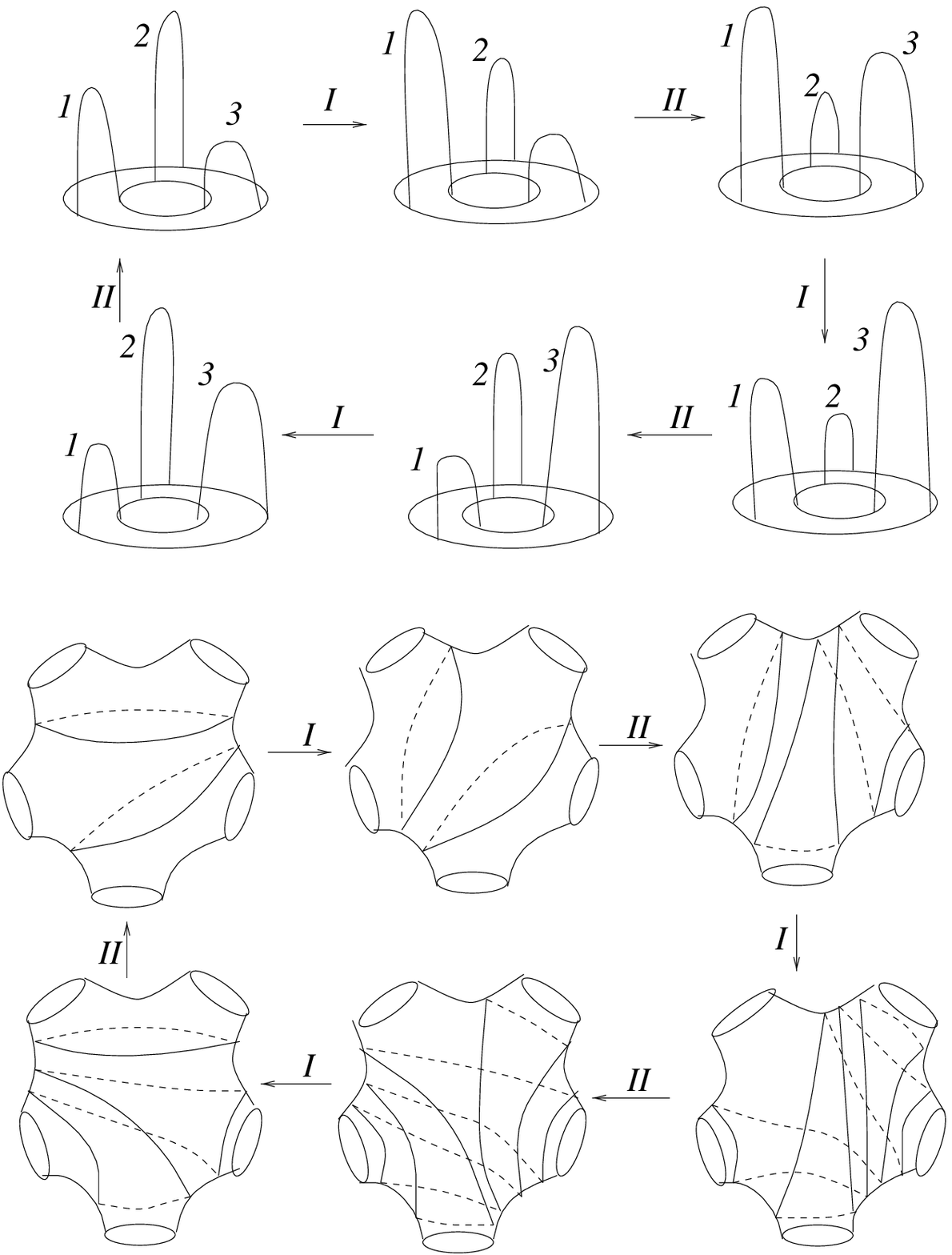}
\caption{Relation 10}
\label{cell10}
\end{figure}

\subsection{Reduction to  fundamental moves and relations}

\begin{figure}[htbp]
\centering
\leavevmode
\epsfxsize=4in
\epsfysize=2.2in
\epsfbox{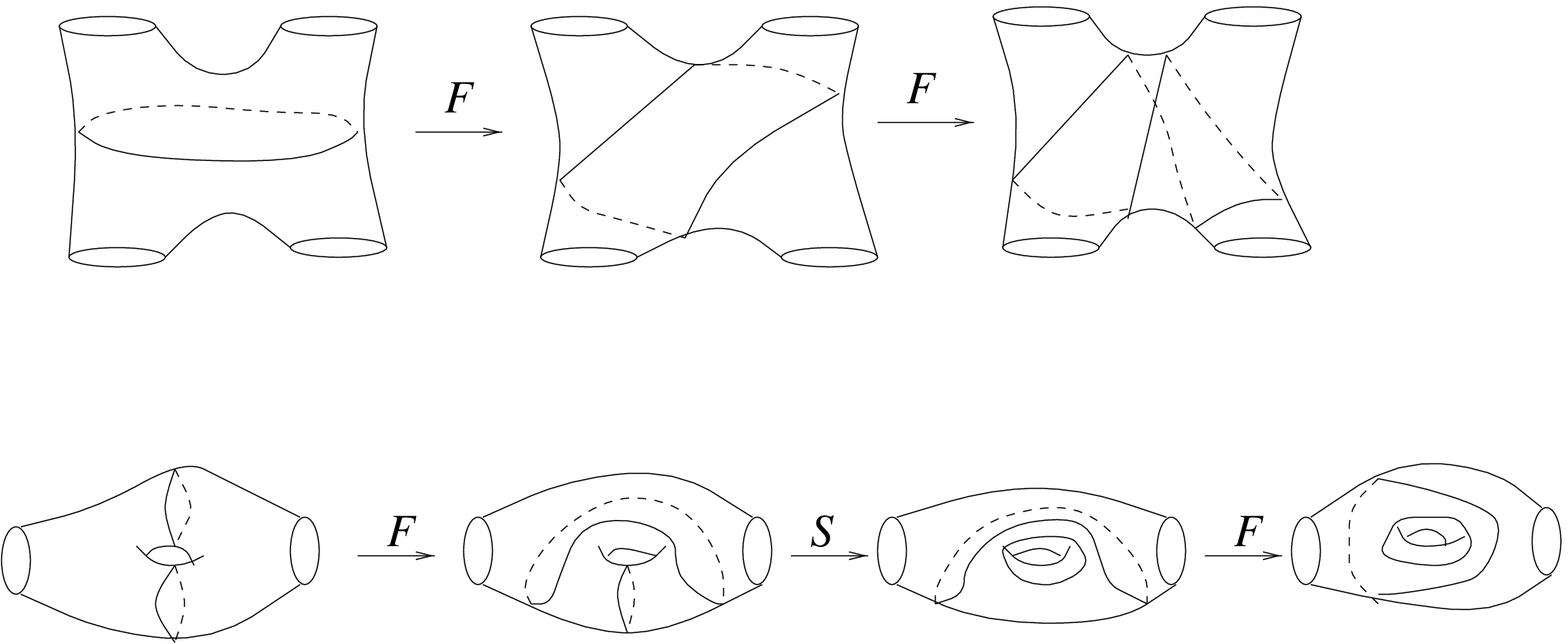}
\caption{Reduction of II and III to $F$ and $S$}
\label{redfs}
\end{figure}

As it is customary in topological quantum field theory, we will
denote the move $I$ by $F$ and the move $IV$ by $S$. The other two
moves can be reduced to these two as seen in 
 Fig.~\ref{redfs}.

Regarding the relations between the moves, recall that Moore and Seiberg
\cite{MoSe} predicted a much smaller number of equations. The reduction
to these is the content of the following proposition.

\begin{figure}[htbp]
\centering
\leavevmode
\epsfxsize=5.5in
\epsfysize=4in
\epsfbox{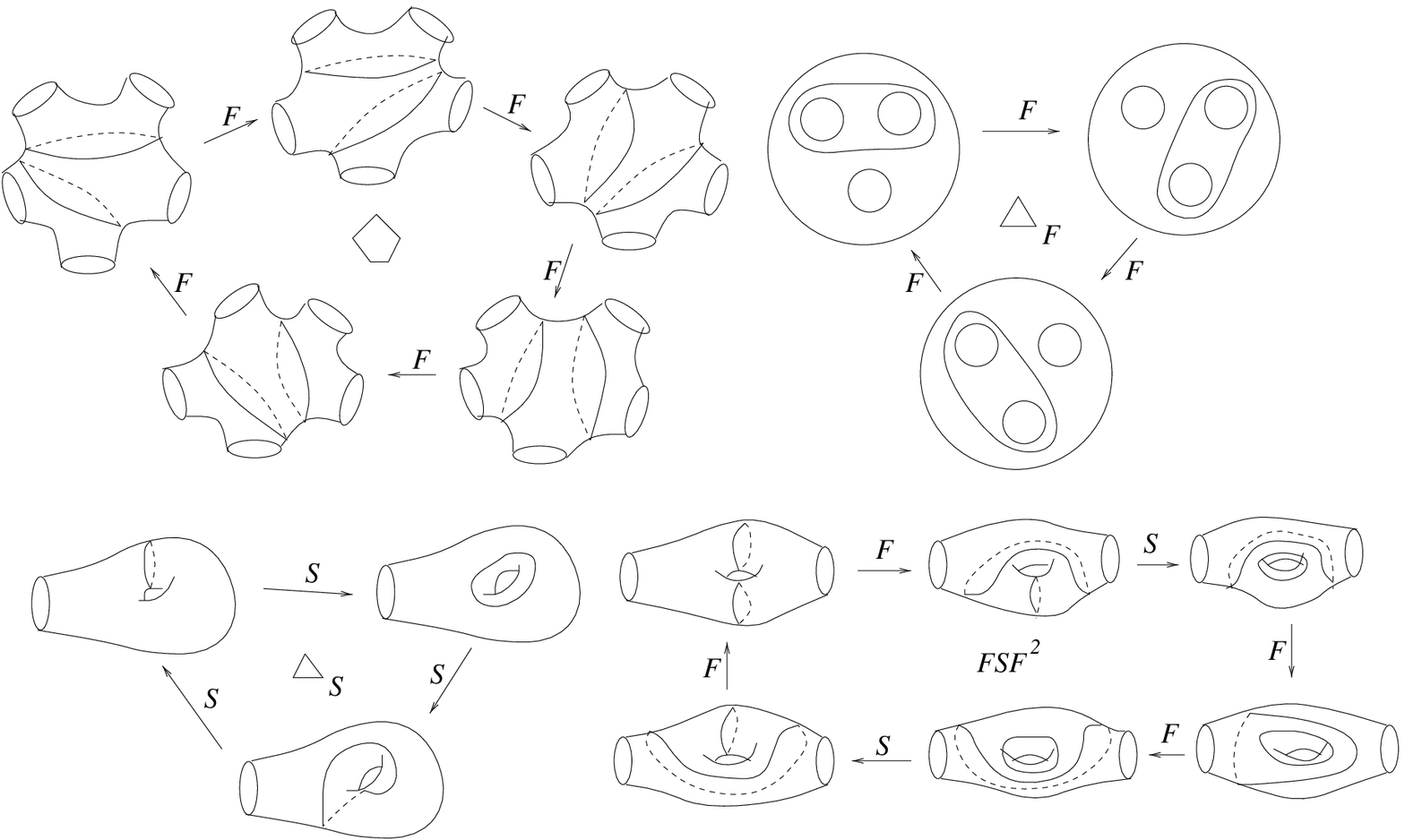}
\caption{Fundamental cells}
\label{fundcell}
\end{figure}

\newtheorem{prop0}[proposition]{Proposition}
\begin{prop0}
Each of the cells arising from triangle singularities can
be decomposed into some of  the four fundamental cells described
in Fig.~\ref{fundcell} and the commutativity 
DC-cells. 
\end{prop0}

Some remarks before we proceed with the proof.
 We have to show that each of the 
seven cells from above decomposes as a union of fundamental cells. 
As stated, this is not quite true, since we must  add some
other cells, the DC-cells, which express the disjoint commutativity 
between $F$ and $S$. Roughly speaking two operations (like $F$ and
$S$) with disjoint supports  commute with each other. 
The squares expressing the commutation are the  DC-cells.  
The reason we need to consider these DC-cells is the 
fact that we didn't take into account the tensorial structure for the 
moment and are thereby working  with a fixed surface.

{\em Proof:} The decompositions are presented in Fig. 20 
through Fig.~\ref{decomp10}. 

\begin{figure}[htbp]
\centering
\leavevmode
\epsfxsize=3.5in
\epsfysize=1.6in
\epsfbox{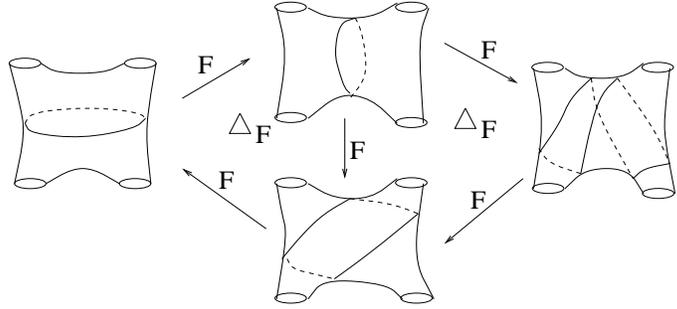}
\caption{Decomposition of $F^2$}
\label{fsquare}
\end{figure}

Let us point out that, as shown in Fig.~\ref{fsquare}, the possibility
of decomposing the cells into fundamental cells does not depend on the
way we expand the moves $II, III$ (the cells in
Fig.~\ref{fsquare} provide a homotopy between the two ways to expand
the moves).

\begin{figure}[htbp]
\centering
\leavevmode
\epsfxsize=4in
\epsfysize=3in
\epsfbox{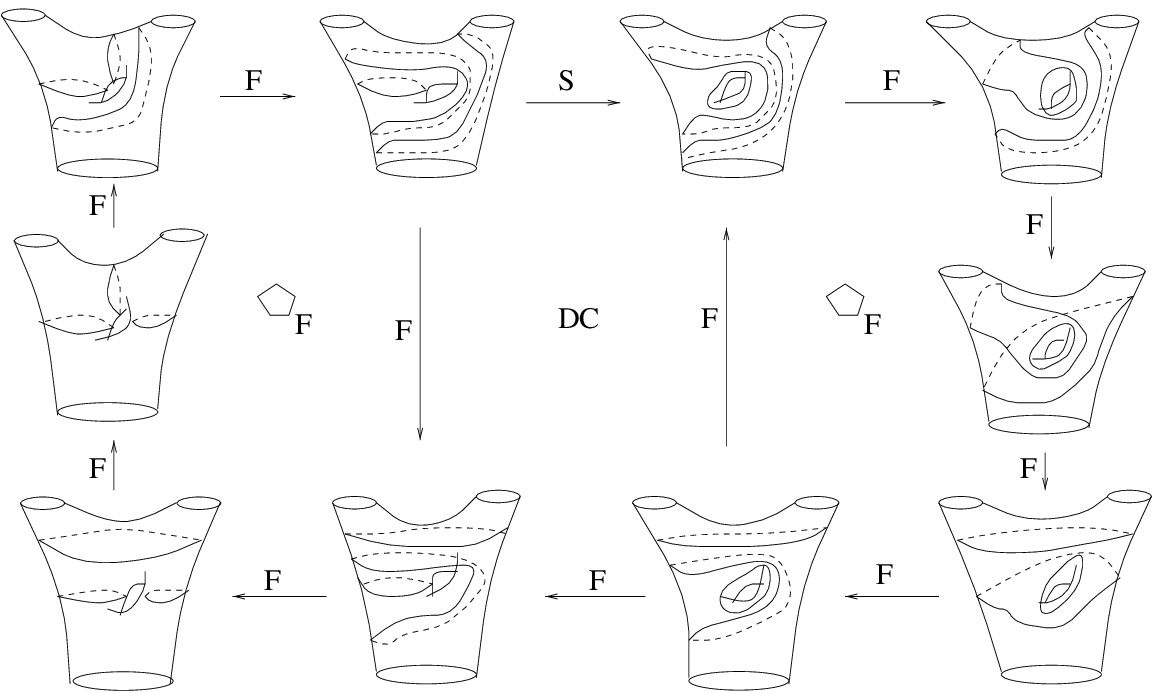}
\label{decomp2}
\caption{Decomposition of cell 2}
\end{figure}

\begin{figure}[htbp]
\centering
\leavevmode
\epsfxsize=4.5in
\epsfysize=4.5in
\epsfbox{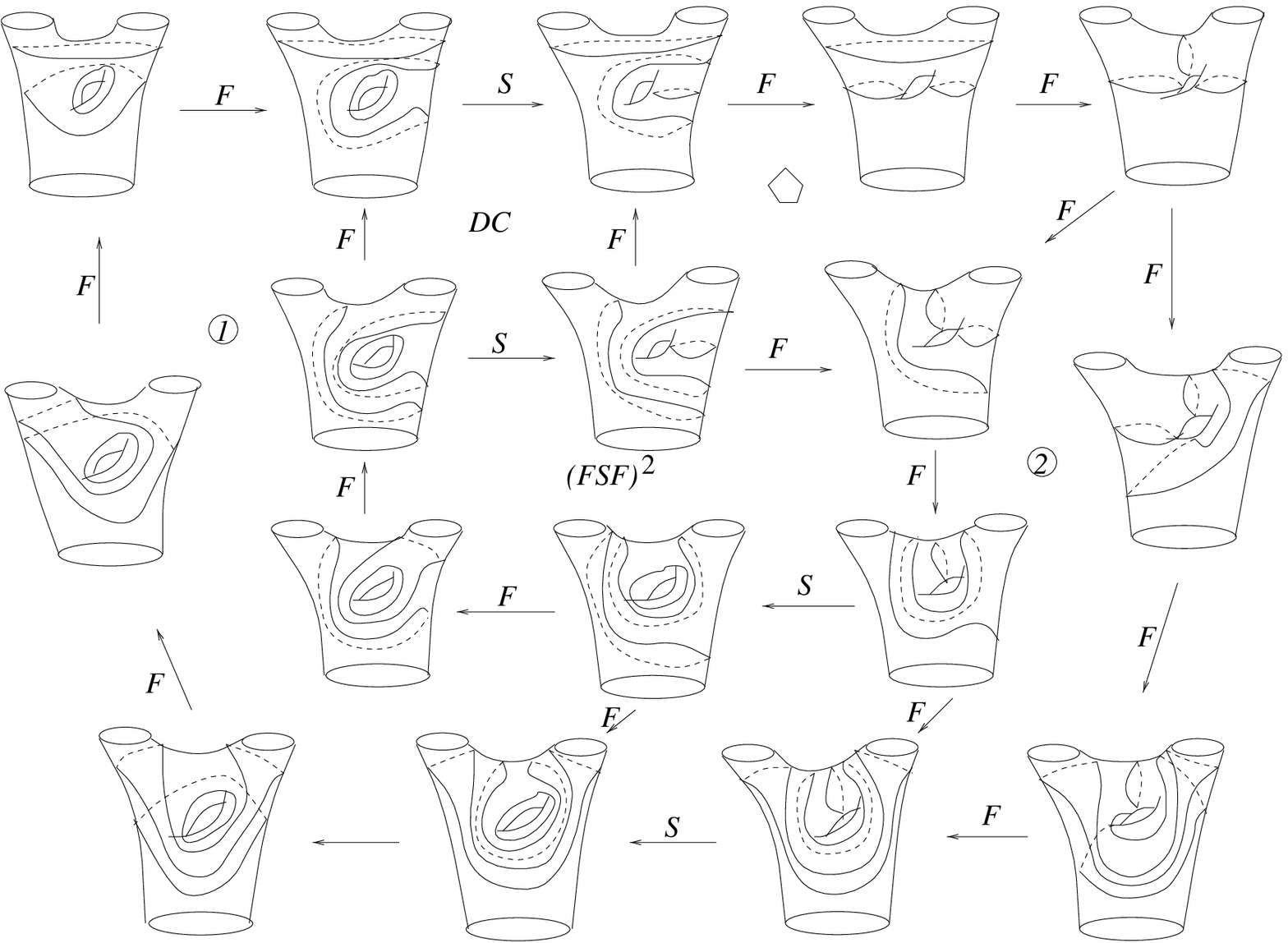}
\caption{Decomposition of cell 3}

\centering
\leavevmode
\epsfxsize=3.5in
\epsfysize=3.5in
\epsfbox{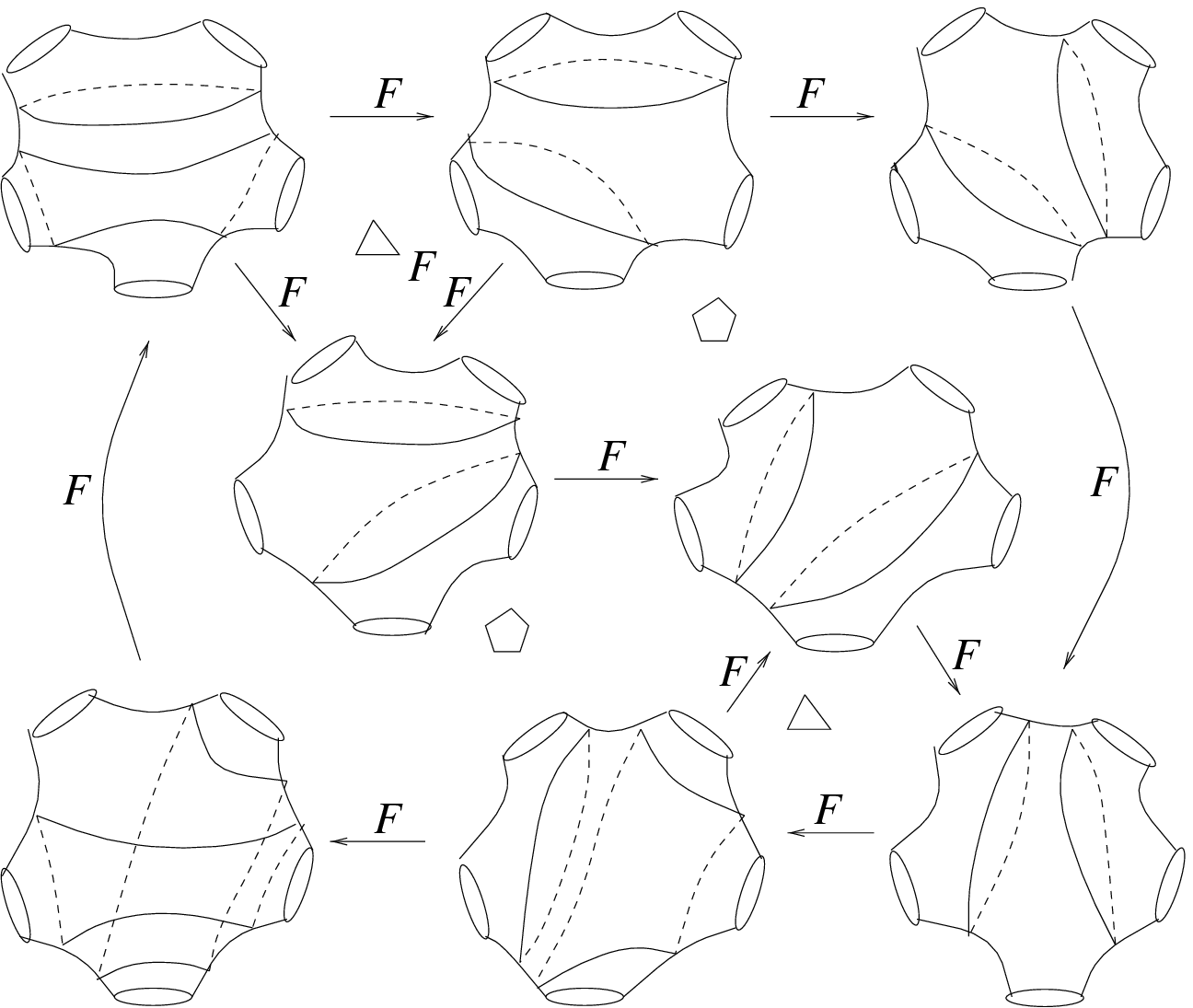}
\caption{Auxiliary cell for decomposition of cell 3}
\end{figure}

\begin{figure}[htbp]
\centering
\leavevmode
\epsfxsize=5in
\epsfysize=5in
\epsfbox{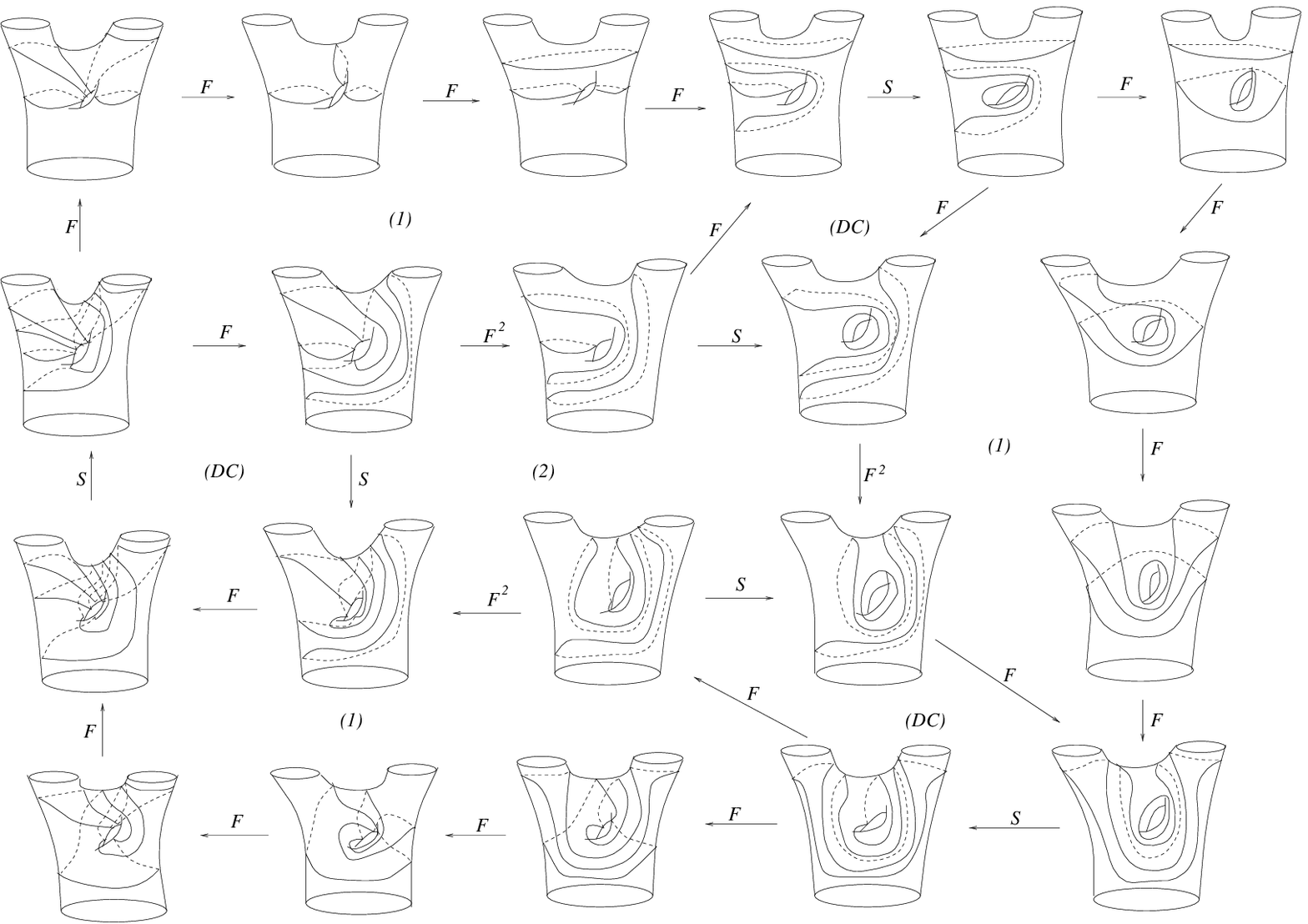}
\caption{Decomposition of cell 4}
\end{figure}

\begin{figure}[htbp]
\centering
\leavevmode
\epsfxsize=3.5in
\epsfysize=3.5in
\epsfbox{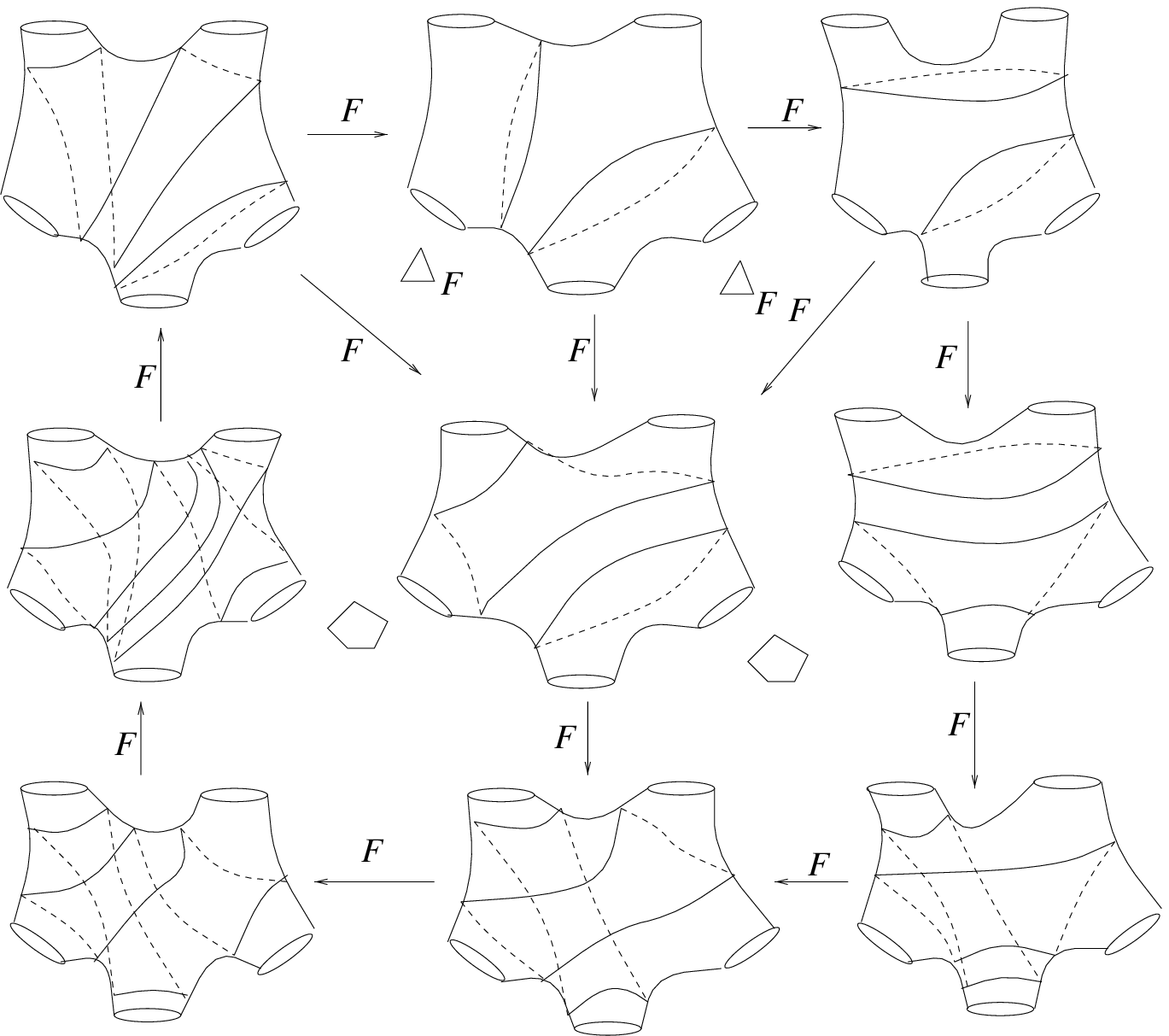}
\caption{Auxiliary cell for decomposition of cell 4}

\centering
\leavevmode
\epsfxsize=4in
\epsfysize=1.8in
\epsfbox{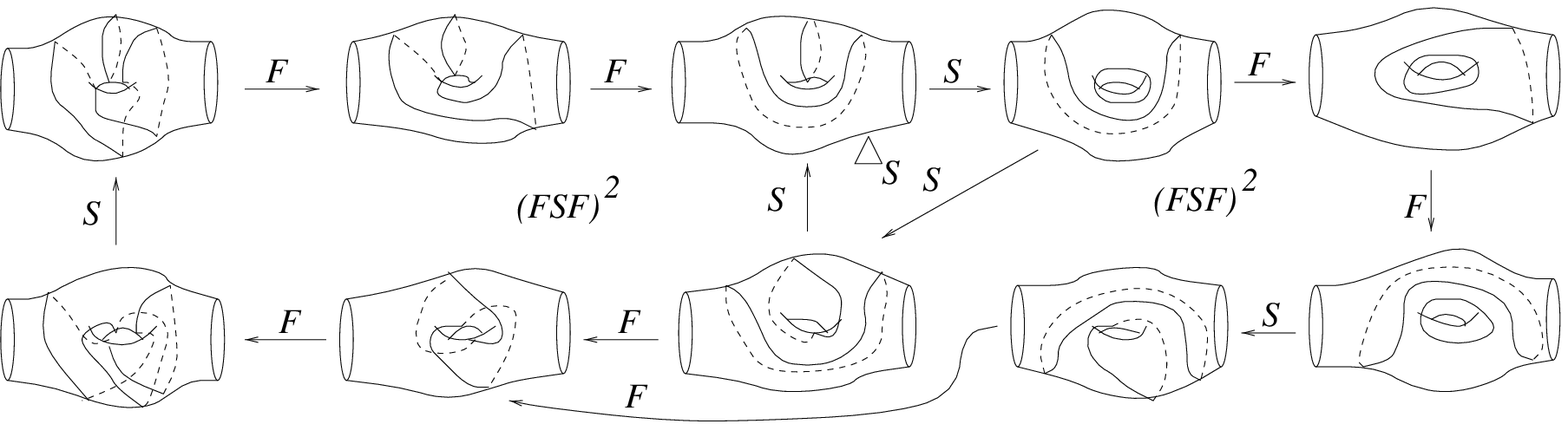}
\caption{Auxiliary cell for decomposition of cell 4}
\end{figure}

\begin{figure}[htbp]
\centering
\leavevmode
\epsfxsize=4in
\epsfysize=3in
\epsfbox{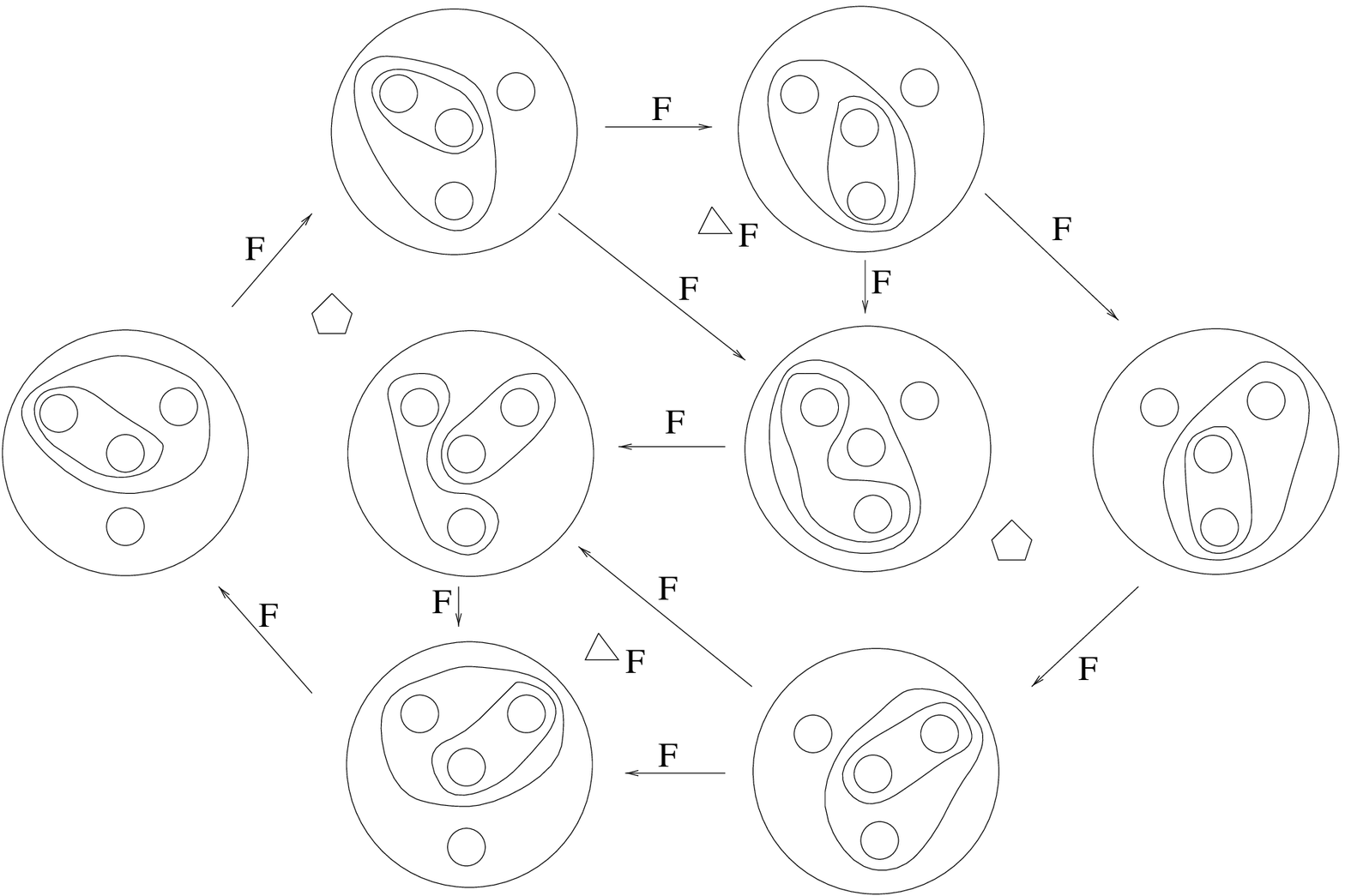}
\caption{Decomposition of cell 5}

\centering
\leavevmode
\epsfxsize=3.5in
\epsfysize=3.5in
\epsfbox{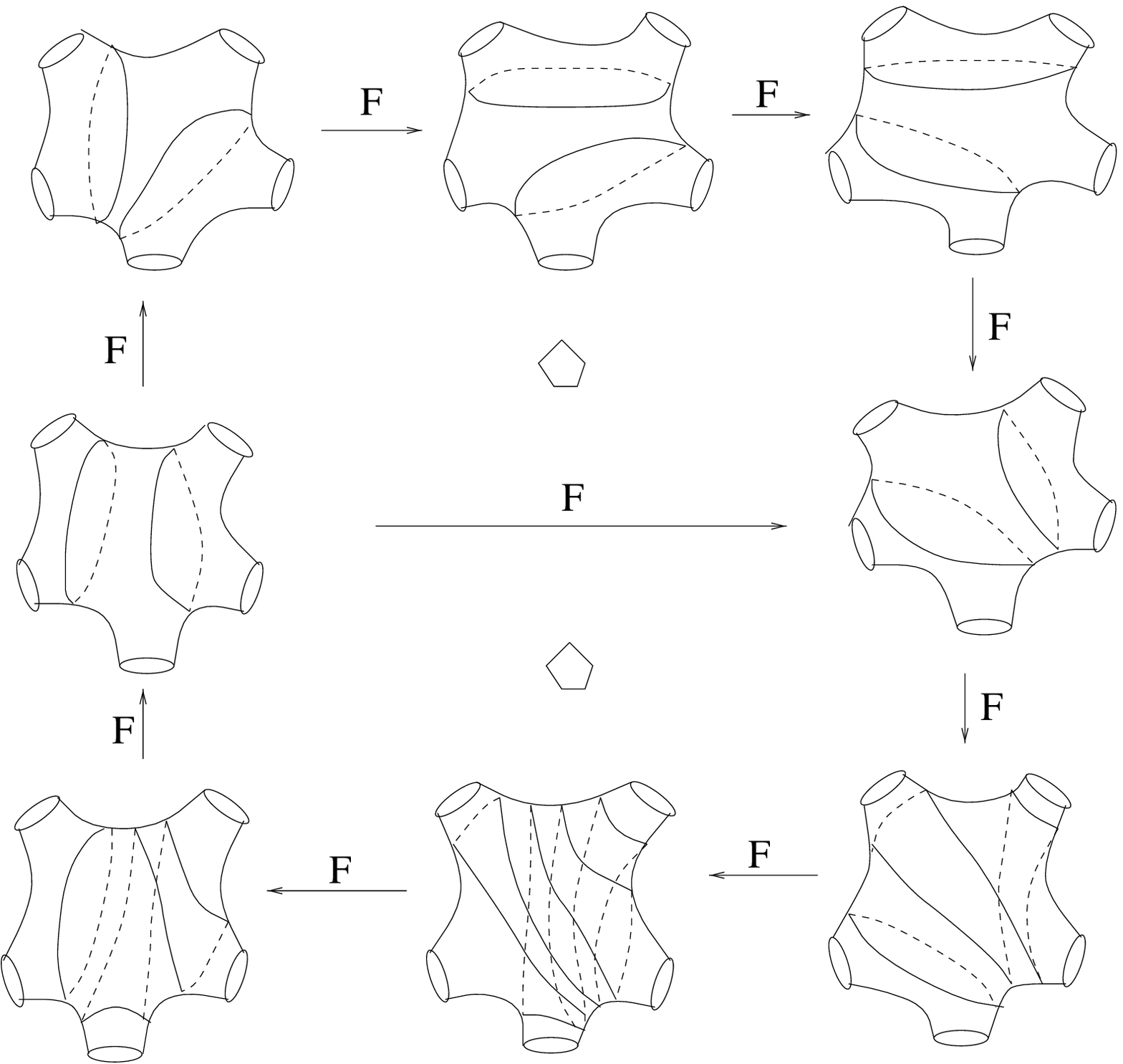}
\caption{Decomposition of cell 9}
\end{figure}

\begin{figure}[htbp]
\centering
\leavevmode
\epsfxsize=4in
\epsfysize=4in
\epsfbox{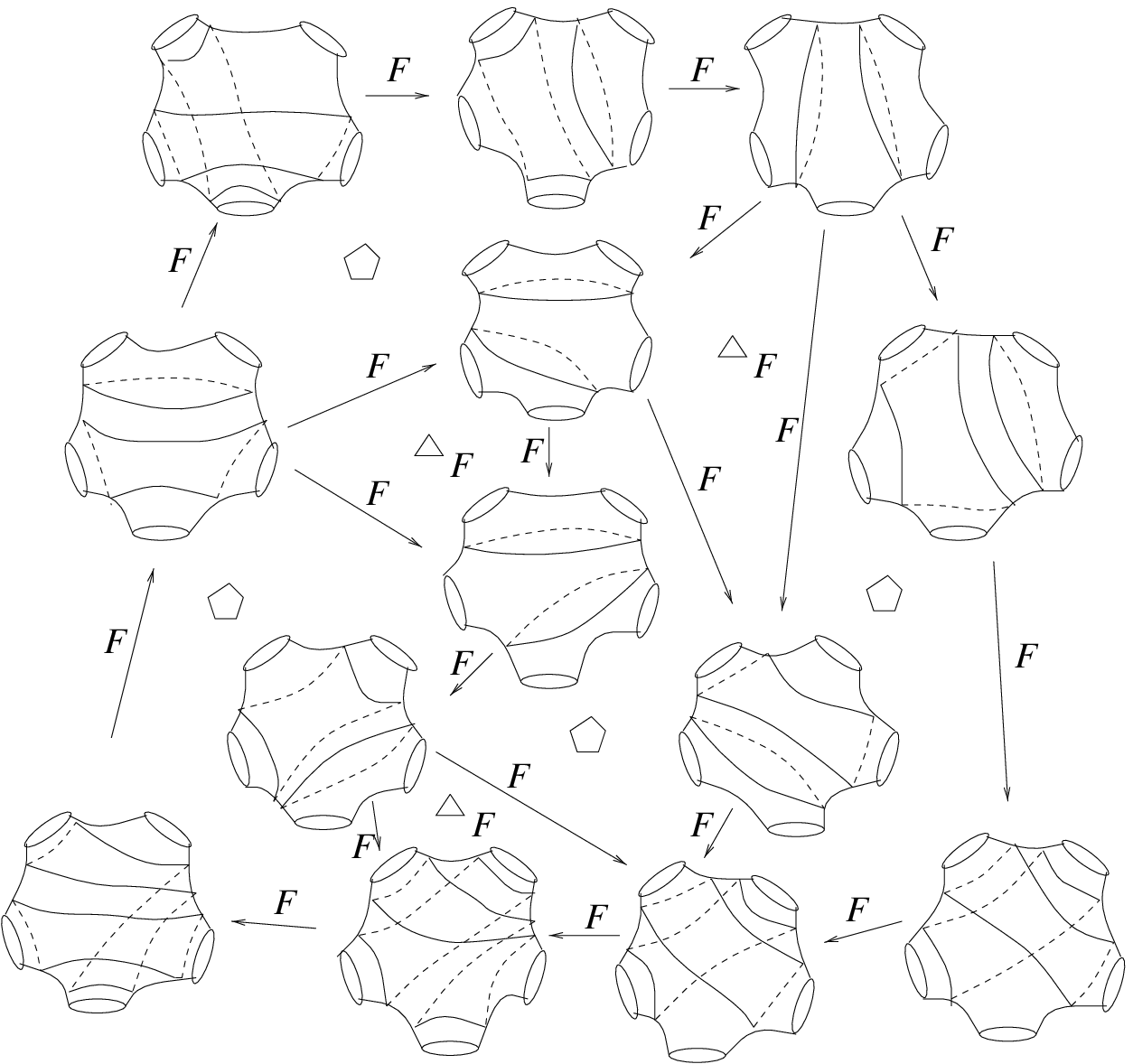}
\caption{Decomposition of cell 10}
\label{decomp10}
\end{figure}

Let us consider now a CW-complex of dimension 2, encoding all 
the informations about markings and transformations. 
The vertices of the complex $\Gamma_0(\Sigma)$ are all 
 possible  markings on the surface  $\Sigma$ (with or without boundary). 
There is an edge (unambiguously denoted  $F$ or $S$) 
between two vertices if the corresponding transformation 
$F$ (respectively $S$) relates the respective markings. 
To this complex we attach the 2-cells described in 
Fig.~\ref{fundcell} and the DC-cells.
 
\newtheorem{prop2}[proposition]{Proposition}
\begin{prop2} 
The complex $\Gamma_0(\Sigma)$ is connected and 
simply-connected. 
\end{prop2}

{\em Proof:}
The result is a consequence of
Proposition 1.3 and Proposition 4.2. More precisely, from Proposition 3.1
we get that $\tilde{\Gamma}_0(\Sigma)$ is connected and simply connected.
 The connectedness of $\Gamma_0(\Sigma)$ follows since each 
move of type $ II$ or $III$ is a composition of 
$F$ and $S$. Notice that the decompositions are not unique.
However, two different decompositions can be homotoped one 
into the other via $F$-triangles and hexagonal $(FSF)^2$-cells respectively
(see Fig.~\ref{fsquare}).  
Also note that the DC-cells made from arbitrary moves  
$I, II, III, IV$ decompose into  DC-cells for $F$ and $S$ according to the 
move decomposition. 
Furthermore the 2-cells in $\tilde{\Gamma}_0(\Sigma)$ are
replaced by their counterparts from $\Gamma_0(\Sigma)$, when the 
respective edges in $\tilde{\Gamma}_0(\Sigma)$ are decomposed. The
latter decompose in $\Gamma_0(\Sigma)$ as unions of fundamental cells 
and DC-cells. This proves the proposition. $\Box$ 
%Remark also that an edge in  $\tilde{\Gamma}_0(\Sigma)$
%which is a composition of $F$ and $S$ has all its decompositions 
%path homotopic when viewed in $\Gamma_0(\Sigma)$. 
%For elementary moves this is clear: a hexagonal cell 
%$(FSF)^2$ gives the homotopy for the two decompositions of  $III$, and 
%two triangles (see the picture) give the homotopy for $II$.   

\subsection{Overmarkings}

Following Walker \cite{Wa} we  call a finite collection
of disjoint simple closed curves lying in the interior
of a surface an overmarking. Such a family of curves
decomposes the surface into disks, annuli and pairs of
pants, decomposition which is also called a DAP-decomposition.

Given a fixed surface, we want to exhibit a 
set of generators and relations for the groupoid of  
transformations of overmarkings. 
A  decomposition
containing  only disks and pairs of pants
is determined  by the level sets of a Morse function.
The disks are semi-local models of points of index $0$ and $2$,
and the pairs of pants are semi-local models of
points of index $1$.
By adding annuli one adds circles that are isotopic to the given circles.

\begin{figure}[htbp]
\centering
\leavevmode
\epsfxsize=4in
\epsfysize=1in
\epsfbox{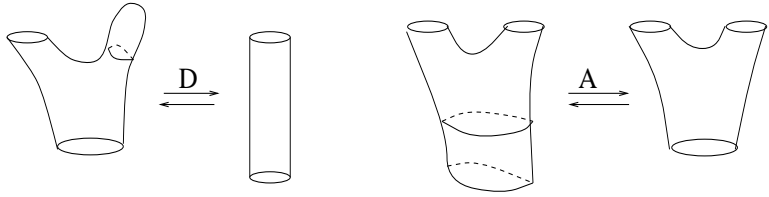}
\caption{Birth-death move}
\label{birth}
\end{figure}

Like before, two decompositions can be transformed one into the other 
along a good path, hence the elementary moves come from crossings
of critical points, and by introducing (expanding) or eliminating
(contracting)  a finite number of annuli. 
In addition to the moves described in the previous section,
one has the  moves described in Fig.~\ref{birth}, where we note that the first
comes from a birth or death point.

\begin{figure}[htbp]
\centering
\leavevmode
\epsfxsize=4in
\epsfysize=4in
\epsfbox{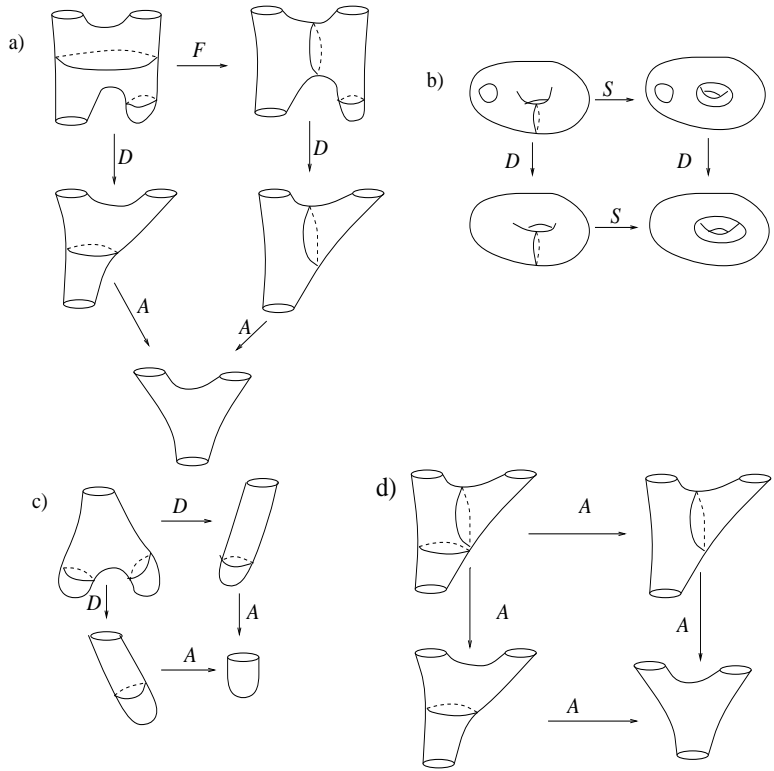}
\caption{Birth-death cells}
\label{bdcell}
\end{figure}

The new 2-cells are the ones produced by birth-death singularities
 (Fig.~\ref{bdcell}. a), b)),
swallow tail singularity (Fig.~\ref{bdcell}. c)), 
and disjoint commutativity.

\begin{figure}
\centering
\leavevmode
\epsfxsize=3in
\epsfysize=3in
\epsfbox{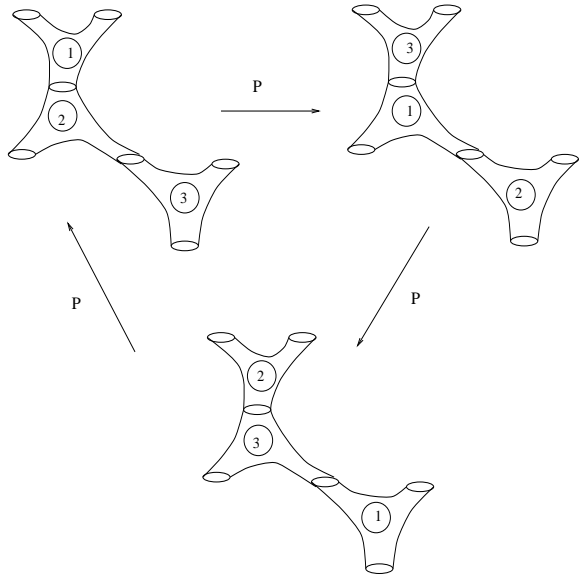}
\caption{An example of a permutation cell}
\label{Pcell}
\end{figure}

Consider now the groupoid of overmarkings and its associated 
2-complex $\Gamma_1(\Sigma)$. Remark that the groupoid  is defined for any
surface $\Sigma$, without restrictions. 
Here the vertices are the overmarkings, the edges 
correspond to the moves $F, S, D, A$ between two interrelated
overmarkings and the 2-cells are four fundamental cells 
from the previous section, together with  those from figure ~\ref{bdcell} 
and all the DC-cells made out of the four elementary moves. 
As a consequence of  the above discussion and Proposition 4.3 we get

\newtheorem{prop3}[proposition]{Proposition}
\begin{prop3}
The 2-complex $\Gamma_1(\Sigma)$ is connected and simply connected. 
\end{prop3}

\begin{figure}
\centering
\leavevmode
\epsfxsize=5in
\epsfysize=5in
\epsfbox{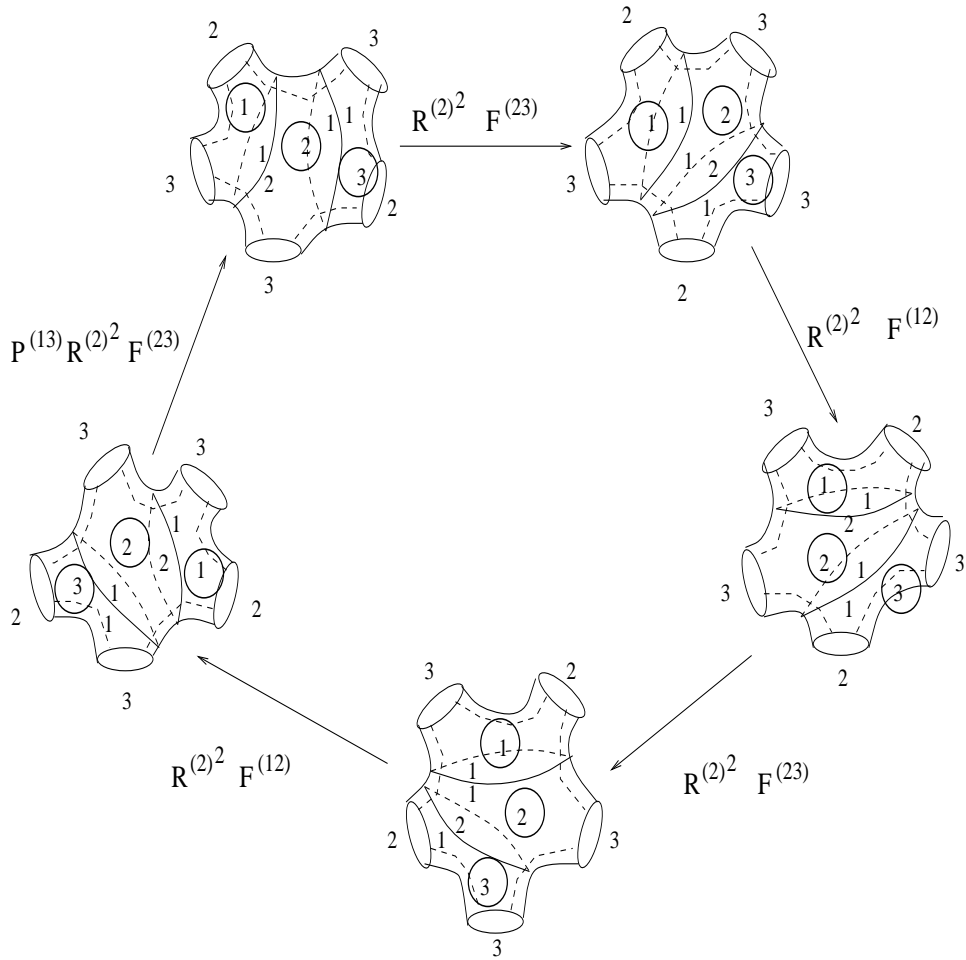}
\caption{The pentagon}
\label{Pentagon}
\end{figure}

\begin{figure}
\centering
\leavevmode
\epsfxsize=5in
\epsfysize=5in
\epsfbox{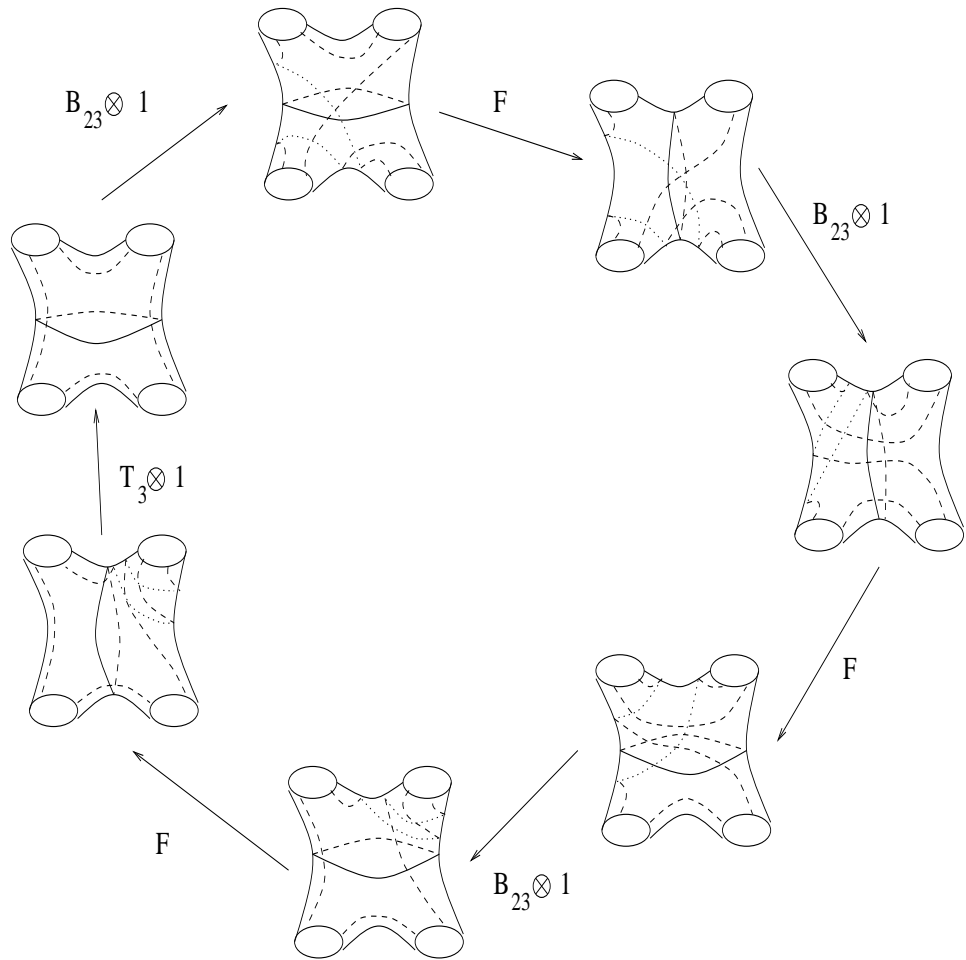}
\caption{The $F$-triangle}
\label{FT}
\end{figure}

\begin{figure}
\centering
\leavevmode
\epsfxsize=5in
\epsfysize=5in
\epsfbox{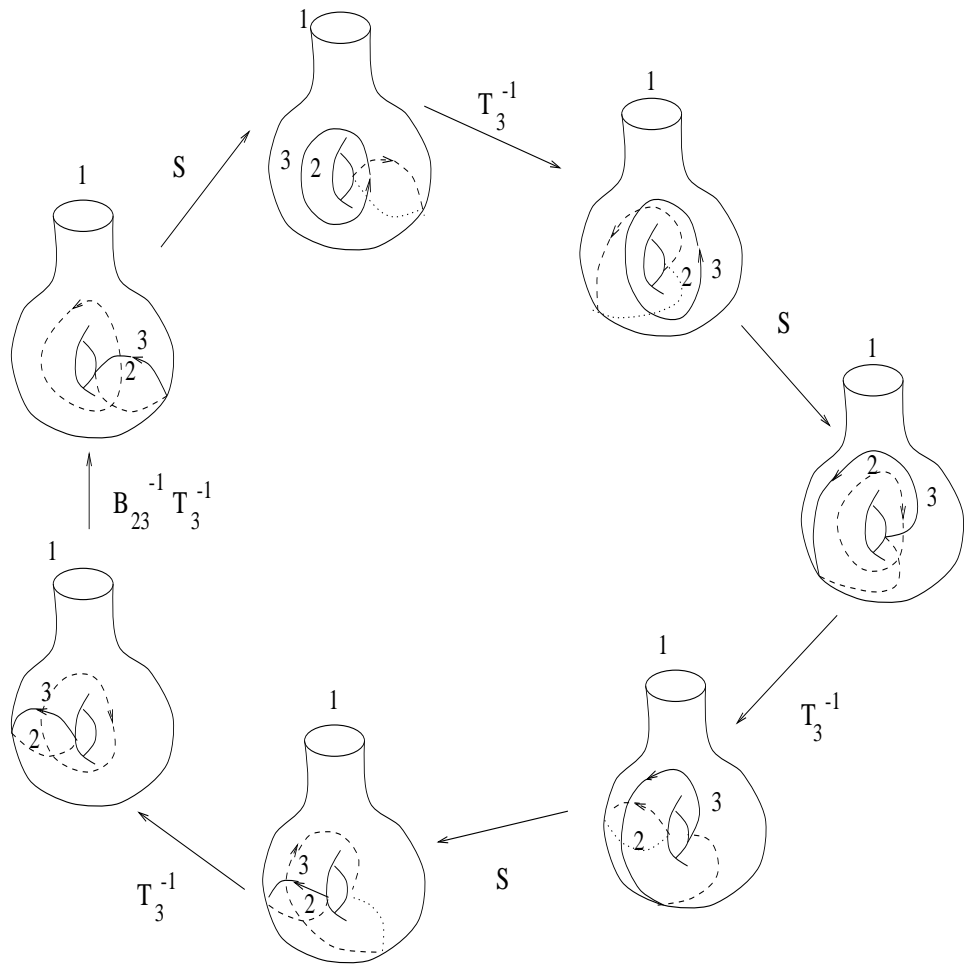}
\caption{The $S$-triangle}
\label{ST}
\end{figure}
\subsection{Rigid structures}

%Remark that upon now the 2-complexes associated were finite, but the
%rigid structures  are infinitely many. This is due to the extra
%structure contained in the seams. This makes the mapping class group 
%act freely on the rigid structures. In some sense this gives a
%construction for a classifying space for the mapping class groups. 

Let us proceed with the proof of Theorem 3.5. 
Fix a surface $\Sigma$ and consider the 2-complex $\Gamma(\Sigma )$
defined in Section 3.
Recall that a rigid structure consists of the following data:
\begin{enumerate}
\item An overmarking $\alpha$ inducing a DAP-decomposition of $\Sigma$.
\item Seams on the elementary surfaces  of the DAP-decomposition.
\item Numberings of the boundary components of these elementary surfaces. 
\item An ordering (segregated according to the topological type)
of the surfaces in the DAP-decomposition.
\end{enumerate}

\begin{figure}
\centering
\leavevmode
\epsfxsize=5in
\epsfysize=6in
\epsfbox{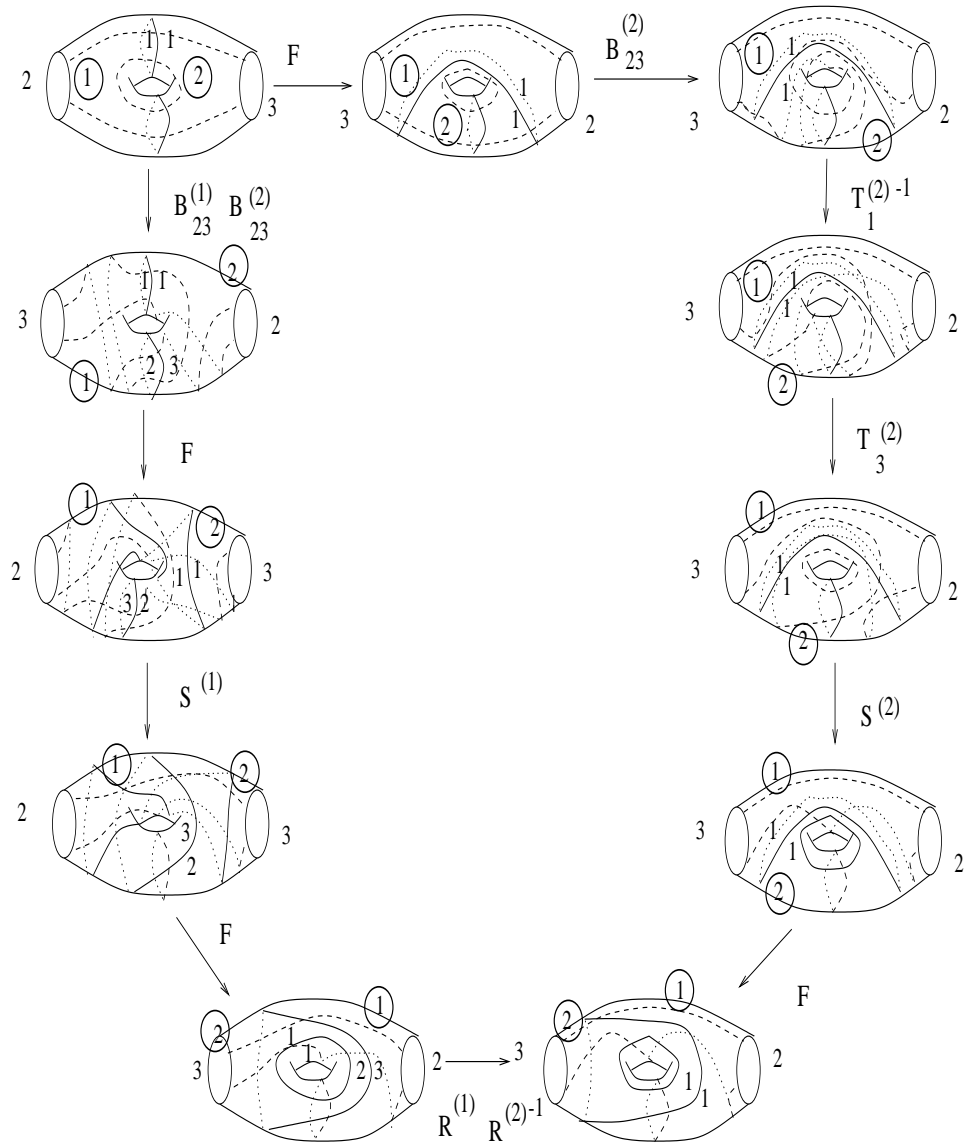}
\caption{The $(FSF)^2$-cell}
\label{FSF}
\end{figure}

The 2-cells of $\Gamma (\Sigma)$ consist  of:
\begin{enumerate}
\item One cell for each cycle of moves of type $P$ when there is 
a corresponding cell in the group of permutations (see Fig.~\ref{Pcell}). 
\item A cell for each relation 1.a)-1.d).  
from Theorem
  3.1. 
\item A cell for relations 2.a) and 2.b) in Theorem 3.1, defining 
inverses. 
\item For each fundamental cell in the complex of
  overmarkings one lifting of this cell at the level of
rigid structures. This means that we consider some labeling of 
one vertex and one system of seams and then keep tracking 
the labeling and the seams all over the 
boundary cell, possibly using the operators $P$ and $R$ which 
permute the numberings, and the twisting operators
 to change the seams configuration. 
\item The DC-cells. 
\item Four  cells, each of which represents a lifting of one of 
the cells made up from $D$ and $A$. This  means that we add 
 seams and numberings 
to the cells from Fig~\ref{bdcell}.  
\end{enumerate}

The liftings of the four fundamental cells are 
shown in Fig.~\ref{Pentagon} through  Fig.~\ref{FSF}.
Note  that in  Fig.~\ref{ST} we have drawn only the closed seam.
The other seam is completely determined up to a twist around
the boundary component labeled by $1$, thus can be ignored.

% that the two seams on the punctured torus are homologous to each
%other, so if we know one of them we can find
%the other. Thus it suffices to draw the closed
%seam with an arrow on it in order to specify the sign of its 
%homology class, to be able to keep track of the moves. We orient
%accordingly the curve to cut open the pair of pants (giving the 
%simplest DAP-decomposition of the holed torus), such that the
%move $S$ have the matrix $\left [\begin{array}{cc}
%0 & -1 \\
%1 & 0 
%\end{array}
%\right ]$. 
%Then the cell lifting the $S$ triangle has the following form:

%We omit in the next picture the numberings automatically determined by
%the existent labels, to save space. One curve from the 
%cut system is doubly labeled by $1$ and the other by $2$ and $3$
%respectively. Then the cell lifting $(FSF)^2$ has the following shape:

To conclude the proof of the theorem consider the 
 canonical map $f:\Gamma(\Sigma)\longrightarrow \Gamma_1(\Sigma)$ 
which forgets about the seams and numberings is cellular and has the
following properties: 
\begin{enumerate}
\item $f^{-1}(z)$ is connected and simply connected for any 
0-cell or 1-cell $z$ of $\Gamma_1(\Sigma)$, 
\item for any 2-cell $y$ of $\Gamma_1(\Sigma)$ there exists 
a 2-cell $x$ in  $\Gamma(\Sigma)$ such that $f(x)=y$. 
\end{enumerate}
Since $\Gamma_1(\Sigma)$ is connected and simply connected, 
standard results in algebraic topology (see also \cite{BK,Wa}) imply that 
$\Gamma(\Sigma)$ is connected and simply connected as well. 

The action of the mapping class group is given by
$f(\Sigma, r)=(\Sigma , f(r))$, where $f(r)$ is the image
of the rigid structure $r$ through the homeomorphism
$f$. Since a homeomorphism is determined up to
isotopy by the image of the rigid structure, the
action is free.
This 
ends the proof of  Theorem 3.5. $\Box$

\section{TQFT and representations of the duality groupoid}
\subsection{Three-dimensional TQFT's}
For the sake of completeness we include some basic definitions 
concerning topological quantum field theories,
 and refer to \cite{Tu} for an extensive treatment. 
Our presentation follows the lines developed in \cite{Fun}, and for simplicity
we 
skip the case of TQFT's with anomaly, which are defined for 3-manifolds 
with an additional structure ($p_1$-structure in \cite{BHMV2} or a 
2-framing with Atiyah's terminology). The latter are related to 
a ${\bf Z}$-extension of our duality groupoid, which  corresponds to the 
central extensions from \cite{MR}.  It is worth  mentioning 
that, although the case of the extended groupoid is analogous to the case of
the non-extended one,
all ``interesting'' TQFT's are extended, 
i.e. they arise   for manifolds with additional 
structure. The situation is entirely similar to the description of 
highest weight representations of $Diff(S^1)$: there are no
highest weight representations of $Vect(S^1)$ but there exist
interesting representations (e. g. Verma modules) for the unique
central extension, namely the Virasoro algebra. This makes the
presence of a central charge necessary. We think that the same
phenomenon holds for the duality groupoid: if one asks the theory 
to have a  unique 
cyclic vacuum vector  (corresponding to the cyclic vector of 
a Verma module),  
and one also requires the theory to be 
 unitary (i.e. to have positive energy) then we must  consider 
a somewhat canonical central extension, which gives  rise to a
$p_1$-structure for 3-manifolds (which is the analogue of the central
charge of the Virasoro algebra). The philosophy behind this
correspondence is a principle known by  the conformal field theory
community, which basically says that {\it representations of the 
tower of mapping class groups (meaning unitary finite dimensional
representations of the duality groupoid giving rise to a TQFT) correspond
to representations of the Virasoro algebra}. 
If one implication is more or less understood, since physicists 
constructed CFT associated to all highest weight 
(positive energy)  irreducible representations of the 
Virasoro algebra  in both the discrete and  the continuous series, 
and thus derived TQFTs via the monodromy of  conformal blocks, 
the other implication is more difficult. We point out 
the references \cite{ADKP,K,BS} where an action of  the Virasoro algebra is
implicitly  carried by the moduli space of curves
with local parameters  around the punctures. 
Detailed proofs and constructions of the conformal blocks 
coming from the highest weight representations are given in 
\cite{ADKP,BS,TK,TUY}. From this data (usually called CFT) we can
construct the TQFT in 3-dimensions (see for instance \cite{Fun2}). 

\newtheorem{tq}{Definition}[section]
\begin{tq}
A TQFT in dimension 3 is a representation of the category of  oriented
3-dimensional cobordisms  into the category of hermitian vector spaces
$V$. 
\end{tq}

In other words a TQFT is a functor assigning to each 
oriented surface $\Sigma$ a hermitian vector space 
$W(\Sigma)$. Then to each cobordism $M$ between the 
surfaces $\partial_+M$ and $\partial M_-$ one  associates 
a linear map $Z(M):W(\partial_+M)\to W(\partial_-M)$. 
This data is subject to the following conditions:
\begin{enumerate}
\item If $\overline{\Sigma}$ is the surface with the 
orientation reversed then 
\[W(\overline{\Sigma})= W({\Sigma})^*.\]
\item 
If $\cup$ denotes here  disjoint union, then the following quantum rule
holds:
\[ W(\Sigma_1\cup \Sigma_2)=W(\Sigma_1)\otimes W(\Sigma_2). \]
\item If the cobordism $M\circ N$ is the composition of the 
cobordisms $M$ and $N$ then 
\[ Z(M\circ N)=Z(M)\circ Z(N). \]
\item We assume that the ground field of the theory is 
${\bf C}$, and thus we put $W(\emptyset)={\bf C}$.    
The theory is called {\em reduced} if $W(S^2)={\bf C}$ holds, and we will 
restrict ourselves to reduced theories in the sequel. 
\item We ask the theory to be topological, which means that 
$Z(M)$ and  $W(\Sigma)$ depend only on the  topological type of the 
manifolds.
\end{enumerate}

The spaces $W_g=W(\Sigma_g)$ associated to a surface of genus $g$ are
also called {\em conformal blocks} in genus $g$. 
The {\em monodromy} of the theory is the series of mapping class 
group representations defined  as follows. 
Assume that $\Sigma_g$ is a fixed standard surface of genus $g$. 
For any $\varphi\in{\cal M}_g$  consider the mapping cylinder 
$C(\varphi)$,  and set 
\[ \rho_g(\varphi)= Z(C(\varphi)) \in End(W_g). \]
The theory is {\em finite dimensional} if all conformal blocks are 
finite dimensional. 
Also the theory is said to be {\em cyclic} if for each $g$ there is one 
orbit of the mapping class group ${\cal M}_gx_g$ which spans linearly 
 $W_g$. It is easy to show that in this case  we can  take $x_g$ to be 
equal to the 
vector $w_g=Z(H_g, id)\subset W(\Sigma_g)$ associated by the TQFT to the 
standard handlebody $H_g$ (the identification of its boundary is 
by the identity map). This vector is called the {\em vacuum} vector in
genus $g$. 

It is shown in \cite{Fun} that any topological invariant $I$ for
closed 3-manifolds defines a series of representations of the 
mapping class group which extends canonically to a TQFT. This is  the
maximal TQFT associated to the invariant $I$.
An important fact to be mentioned is that
the  maximal TQFT is always
cyclic. 
 Notice that the
maximal TQFT is 
uniquely defined, but the same invariant for closed manifolds can arise 
from several distinct TQFT's. 
Starting with a certain 
 invariant of closed manifolds, 
which is the restriction of a TQFT, $Z_0$, 
 and using  the   method described above, one  derives another  TQFT, 
$Z$ which is cyclic and contains basically the same topological
information.

 As an example,  the $sl_2({\bf C})$-TQFT described 
by  Kirby and Melvin \cite{KiMe} is 
not cyclic (and therefore not maximal). The BHMV theories 
(\cite{BHMV1,BHMV2}) give rise to the same invariants for closed manifolds, and
are maximal by construction.
 Notice that, in particular,  all TQFT's  which are
not cyclic induce representations
of the mapping class group which are not
irreducible. 

\subsection{Representations of the mapping class group  and TQFT}
We know that any TQFT determines a series of representations of 
${\cal M}_g$. The converse is also true since the  latter determines  
the  TQFT. We assume from now on that the TQFT,  $Z$ is  
cyclic. 

A cyclic TQFT has more structure hidden in the conformal blocks: 
for instance using the connected sum of 3-manifolds along the 
boundary we derive that there is a natural (injective) homomorphism 
\[ W_g\otimes W_h\hookrightarrow W_{g+h},  \]
induced by $w_{g+h}=w_g\otimes w_h$. 
Call  the sequence $v_g\in W_g$ a sequence of {\it vacuum vectors} for
the representations $\rho_g$ if $v_{g+h}=v_g\otimes v_h$ and 
$\rho_{g}({\cal M}_g^+)v_g=v_g$, for any $g$. Here ${\cal M}_g^+$
denotes the mapping class group elements which arise 
from homeomorphisms of $\Sigma_g$ extending over the handlebody $H_g$.
The TQFT  is said to have {\it  an unique vacuum} if the vector $v_1$ 
which is $\rho_{g}({\cal M}_1^+)$ is unique and therefore equal to 
$w_1$ (up to a scalar). 

Let us point out  that the spaces $W_g$ are naturally endowed
with a hermitian form $<,>:W_g\otimes W_g\to {\bf C}$, given by 
\[ <X,Y> = Z(X\cup \overline{Y}). \]
Here $X, Y\in W_g$  are linear combinations of elements 
$(H_g,\varphi)=\rho_g(\varphi)w_g\in W_g$. This follows from
the fact that the theory is cyclic. It suffices then to consider 
the case $X=(H_g,\varphi_1)$  
 and $Y=(H_g,\varphi_2)$, 
where $\varphi_j\in {\cal M}_g$. The right hand side is the invarinat of the
manifold $X\cup \overline{Y}$,  obtained by  gluing two handlebodies  
$H_g\cup{\varphi_1\varphi_2^{-1}}\overline{H_g}$, where
${H_g}$ is the 
standard handlebody of genus $g$.  

The TQFT is {\em non-degenerate} if the hermitian form $<,>$ is 
non-degenerate for all genera $g$. Obviously we can replace 
$W(\Sigma_g)$ by $W_g/\ker <,>$ in order to make the theory
non-degenerate, without really changing its topological content. 
In particular the invariants of closed manifolds are the same in both
TQFT's. 
A TQFT is called {\it unitary} if the non-degenerate 
form $<,>$ is positive definite. In this
case  the mapping class groups act on $W_g$ by 
unitary operators since the hermitian form is ${\cal M}_g$-invariant. 
The unitarity of the main examples of TQFT's, e.g. the  $SU(n)$-TQFT, 
is the key point in many applications, for
instance in obtaining  the lower bounds for the genus of a knot. 

We can now recover the invariant $Z$ associated to closed 3-manifolds 
from the representations $\rho_*$ and the hermitian form. 
The relationship between the representation and the invariant arising
in the SU(2)-TQFT was obtained in \cite{R}. 
Specifically we have the
following result from \cite{Fun}:
\newtheorem{lapte}{Propositon}[section]
\begin{lapte}
Let $M=H_g\cup_{\varphi}\overline{H_g}$ be a Heegaard splitting of the
manifold $M$, where $\varphi$ denotes the gluing homeomorphism. 
Then 
\[ Z(M)=d^{-g} <\rho_g(\varphi)w_g,\overline{w_g}>,  \]
where $d=<\rho_1(\iota)w_1,\overline{w_1}>$ is a normalization factor, 
$\iota$ being the gluing map associated to the standard Heegaard
splitting of $S^3$ into two solid tori (i.e.  we choose 
the standard basis in the homology of the torus and 
$\iota= \left [ \begin{array}{cc}
                0 & -1 \\
                1 & 0 
                \end{array}
         \right ] \in SL_2({\bf Z})$.) 
\end{lapte}

\subsection{Proof of  Theorem 3.2}
 
The proof of the theorem consists basically of the construction of the
representation of the $2$-groupoid.
We  outline first the structure of a cyclic unitary TQFT with 
an unique vacuum, along  the lines of  \cite{Fun}. 

We define the primary conformal blocks $W^{i}_{jk}$ as follows. 
Let $r$ be a rigid structure on the 
surface $\Sigma_g$, made from the  pants decomposition $c$, the seams 
and the various numberings.  On each trinion, the set of three seams 
connecting the boundary components can be uniquely identified 
with the boundary of the neighborhood of the graph $Y$ embedded in the
pants. More precisely, $Y$ is the graph topologically isomorphic to the
letter $Y$ and it is properly embedded in the trinion. Let us then
consider one such graph for each trinion and then their union is a 
3-valent graph $\Gamma$ of genus $g$ (possibly with some additional leaves). 
The graph $\Gamma$ encodes all the informations carried by the set of
seams. Notice that this is naturally embedded in the surface
$\Sigma_g$ and thus there is a natural cyclic order on the edges around
each vertex. In the language of \cite{Fun} we have a rigid structure,
in fact equivalent to those considered in this paper. 
  
The label  set $A$ is the set of eigenvalues and their 
inverses for  of the Dehn twists $T_{c_j}$ around the curves $c_j$  in the 
pants decomposition. Fix a vertex $v$ in the graph $\Gamma$ whose 
adjacent edges are $e_1,e_2,e_3$ which are dual to the 
curves $c_1,c_2,c_3$ bounding a pair of pants $p$. 
Let us consider a vector $w_g(i_1,i_2,i_3)\in W_g$ such that 
$\rho_g(T_{c_{\alpha}})$ has $w_g(i,j,k)$ an eigenvector of 
eigenvalue $i_j$ if $\alpha=j$ and $1$ otherwise. 
The span $W(i_1,i_2,i_3)$ of the orbit of the  
$w_g(i_1,i_2,i_3)$ by those elements of 
$\rho_g({\cal M}_g$ which come from homeomorphisms of ${\cal M}_g$ 
which have the support on  the trinion $p$. 
It was proved in \cite{Fun} that $W(i,j,k)$ does not depend of the 
choice of the vertex $v$, the rigid structure and  the genus of the
surface $\Sigma_g$. Moreover, with the  convention of  adding orientations
to the edges of $\Gamma$ such that on each vertex there are two
incoming and one outcoming edge, one obtains this way a
 well-defined space 
 denoted $W^{i_1}_{i_2i_3}$. Of course one  should add possible leaves
to the graph, labeled all by $1$, which correspond to capping of the 
surface with disks. These  correspond to the moves of type $D$.

We define now  graphical rules of associating vector spaces to
(partially labeled) graphs. Consider an oriented trivalent graph 
 whose edges are labeled. Each internal vertex has two incoming edges and one 
outgoing edge. Consider the counter-clockwise  cyclic order of  the incident
edges of a vertex.  If  we label 
the edges by elements of the set $A$ 
there is a non-ambiguous way to associate to each 
internal vertex a vector space $W^{\nu}_{\lambda\mu}$ such that $\nu$ is the label of 
the outgoing edge, and the triple $(\lambda,\mu,\nu)$ is cyclicly ordered. 
We associate to the whole labeled graph $\Gamma$  the tensor product of all 
 spaces 
associated to vertices. Finally if the graph has some of its edges
with fixed  labels,  
take the sum of all the spaces obtained by the above construction, over all 
possible labelings of the remaining edges and call this space $W(\Gamma)$.
Remark that these conventions make sense for an arbitrary trivalent graph. 

For  a closed (oriented) surface $\Sigma_g$ of genus 
$g$, endowed with the rigid structure $r$ 
we  consider the subjacent  pants decomposition. 
We associate to the rigid structure
 the trivalent graph $\Gamma\subset \Sigma_g$,
 whose regular  neighborhood contains all seams. 
Notice that the rigid structure may contain an overmarking  instead of
a pants decompositions. Then all the circles which bound (equivalently
all the leaves in the graph $\Gamma$) are labeled by the unit $1$. 
In \cite{Fun} it is proved that:
\newtheorem{lapte2}[lapte]{Proposition}
\begin{lapte2}
For a TQFT, $Z$ which is  unitary, cyclic and has unique vacuum the
conformal blocks decompose in terms of the primary conformal blocks
$W^{i}_{jk}$. This means that there exists a set of labels $A$ and a
set of vector spaces $W^{i}_{jk}$ with the property that  
for a  rigid structure $r$ and a choice of a basis 
in each vector space $W^{i}_{jk}$ (i.e. we fix the internal
symmetries)  there exists a canonical 
isomorphisms $\Phi(\sigma):W_g\to W(\Gamma)$. 
\end{lapte2}
Returning to the proof, consider two rigid structures $\sigma$ and $\sigma'$ 
such that $\sigma=\varphi\sigma'$ for some 
$\varphi\in{\cal M}_g$. It follows  that 
the endomorphism of $W_g$ given by $\Phi(\sigma)\Phi(\sigma')^{-1}$
is equal to $\rho_g(\varphi)$. For two arbitrary 
rigid structures $\sigma$ and $\sigma'$, there exists a 
unique element $[\sigma,\sigma']$ of the duality groupoid  sending 
$\sigma$ into $\sigma'$. We associate  to 
this element  $[\sigma,\sigma']$ the isomorphisms of vector spaces 
$\Phi([\sigma,\sigma'])=\Phi(\sigma')^{-1}\Phi(\sigma): W(\Gamma)\to W(\Gamma')$.
 We define in this way a representation of the 
2-groupoid. The only remark to add 
is that  this map is local, that is if $\sigma$, and $\sigma'$ are 
identical out of the rigid structures $\sigma_0$ and respectively 
$\sigma_0'$ on some subsurface then 
$\Phi([\sigma,\sigma'])= \Phi([\sigma_0,\sigma_0'])\otimes 1$, where 
$1$ acts on the other factors of the tensor product, not coming from
the subsurface. This proves the theorem.

As a final remark, the representation of the duality groupoid for the
case of the BHMV topological quantum field theory \cite{BHMV1,BHMV2}
was given in \cite{Gel2}. The fact that the topological 
quantum field theory from \cite{KiMe} is not cyclic produces a
sign obstruction for constructing a representation of the
duality groupoid in this case. One can construct a representation
of an extension of this groupoid obtained by adding auxiliary 
structure on the boundary of $3$-manifolds \cite{Gel}.

\newpage

\bibliographystyle{plain}
     
\end{document}